\documentclass[twocolumn]{autart}

\newtheorem{example}{Example}

\usepackage{graphicx}          
\usepackage{tikzit}
\usepackage{amsmath}
\usepackage{amsfonts}
\usepackage{amssymb}
\usepackage[hidelinks]{hyperref}
\usepackage{float}
\usepackage{subcaption}
\usepackage{easyReview}
\usepackage{enumitem}

\tikzstyle{Axes}=[->, >=stealth, thick]
\tikzstyle{Graph region}=[-, fill={rgb,255: red,191; green,191; blue,191}, draw=none]
\tikzstyle{axes_tick}=[-]
\tikzstyle{distance_arrow}=[<->, >=stealth]
\tikzstyle{anti_region}=[-, fill=white, draw=none]

\usepackage{tikz}
\usepackage[utf8]{inputenc}
\usepackage{pgfplots} 
\usepackage{pgfgantt}
\usepackage{pdflscape}
\pgfplotsset{compat=newest} 
\pgfplotsset{plot coordinates/math parser=false}
\usepackage{tikzscale}
\usepackage{standalone}
\usetikzlibrary{shapes,arrows}
\usetikzlibrary{arrows,calc,positioning,fit}
\journal{X}
\date{11 July 2025}

\usepackage{eso-pic}
\AddToShipoutPictureBG*{%
\AtPageUpperLeft{%
\setlength\unitlength{1in}%
\hspace*{\dimexpr0.5\paperwidth\relax}
\makebox(8,-1.2)[c]{
\begin{tabular}{c c}
T. de Groot \emph{et al.},
A Dissipativity Framework for Constructing Scaled Graphs, \\
uploaded to arXiv  11 July 2025 \\
\end{tabular}}}}
\begin{document}

\begin{frontmatter}

\title{A Dissipativity Framework for Constructing Scaled Graphs} %

\thanks[footnoteinfo]{The research leading to these results has received funding from the European Research Council under the Advanced ERC Grant Agreement PROACTHIS, no. 101055384.}
\author[TUe]{T. de Groot}\ead{t.d.groot2@tue.nl},     
\author[TUe]{W.P.M.H. Heemels}\ead{m.heemels@tue.nl},
\author[TUe]{S.J.A.M. van den Eijnden}\ead{s.j.a.m.v.d.eijnden@tue.nl}               

\address[TUe]{Eindhoven University of Technology, Eindhoven, Netherlands}


\begin{abstract}                           
Scaled relative graphs have been originally introduced in the context of convex optimization and have recently gained attention in the control systems community for the graphical analysis of nonlinear systems. Of particular interest in stability analysis of feedback systems is the scaled graph, a special case of the scaled relative graph. In many ways, scaled graphs can be seen as a generalization of the classical Nyquist plot for linear time-invariant systems, and
facilitate a powerful graphical tool for analyzing nonlinear feedback systems. In their current formulation, however, scaled graphs require characterizing the input-output behaviour of a system for an
uncountable number of inputs. This poses a practical bottleneck in obtaining the scaled graph of a nonlinear system, and currently limits its use. This paper presents a framework grounded in dissipativity for efficiently computing the scaled graph of several important classes of systems, including multivariable linear time-invariant systems, impulsive systems, and piecewise linear systems. The proposed approach leverages novel connections between linear matrix inequalities, integral quadratic constraints, and scaled graphs, and is shown to be exact for specific linear time-invariant systems. The results are accompanied by several examples illustrating the potential and effectiveness of the presented framework.
\end{abstract}
\end{frontmatter}



\section{Introduction}
Scaled relative graphs have been originally introduced in the context of convex optimization, and provide rigorous and intuitive visual proofs for convergence of a wide range of convex optimization algorithms \cite{ryu_scaled_2022}. More recently, the notion of scaled relative graphs has gained traction in the systems and control community \cite{chaffey_scaled_2021}, \cite{chaffey_graphical_2023}. The scaled relative graph provides a representation  of \emph{incremental} input-output properties of a system characterized by gain and phase, and is useful for feedback analysis \cite{chen2025softhardscaledrelative}. In \emph{non-incremental} analysis of feedback systems, the scaled graph, which is a special case of the scaled relative graph, is of particular interest. In many ways, a scaled graph can be seen as a generalization of the classical Nyquist plot for linear-time invariant (LTI) systems, and, therefore, can pave the way to  a powerful graphical framework for the design and analysis of nonlinear feedback systems \cite{chaffey_scaled_2021},\cite{chaffey_graphical_2023}. The application of scaled (relative) graphs has since its introduction been  expanded across various domains in the systems and control literature, including Lur’e systems \cite{krebbekx_srg_2024}, (non)monotone operator theory \cite{chaffey_graphical_2023}, \cite{quan_scaled_2024}, dominance theory \cite{chen2025graphicaldominanceanalysislinear}, model reduction \cite{Chaff22}, multivariable systems \cite{degroot2025exploitingstructuremimoscaled}, \cite{baron-prada_mixed_2025}, \cite{baron-prada_stability_2025}, multi-agent systems \cite{baron2025decentralized}, and reset control systems \cite{van_den_eijnden_scaled_2024}.

Although substantial progress has been made in the theory behind scaled graph analysis, obtaining scaled graphs for nonlinear systems remains a key challenge and currently limits their broader applicability in nonlinear system analysis and controller design. For specific types of LTI systems, scaled graphs can be derived directly from generalized Nyquist diagrams \cite{chaffey_graphical_2023}, \cite{pates_scaled_2021}, but no such characterization exists for nonlinear systems. In essence, the main difficulty in constructing scaled graphs comes from the fact that this requires characterizing the input-output behaviour of the system for an uncountable set of input signals. This issue was recently addressed in \cite{van_den_eijnden_scaled_2024}, where a computational method for computing an over-approximation of the scaled graph of the specific class of reset control systems was proposed. A key ingredient in that approach is the use of the generalized Kalman-Yakubovich-Popov (KYP) lemma \cite{iwasaki_generalized_2005} to determine input-specific input-output properties of reset systems, akin to characterizing ``mixed'' systems, see also \cite{CHAFFEY2022105198}, \cite{griggs2009interconnections}, \cite{EijCha}. Each input-specific property can be verified by solving specific sets of linear matrix inequalities (LMIs), and maps to a region in the complex plane that partly covers the scaled graph of the reset system. Taking the union of these regions then forms an over-approximation of the scaled graph, which leads to sufficient conditions for guaranteeing specific stability or input-output properties when used in a feedback setting. Although effective, the union-based approach in \cite{van_den_eijnden_scaled_2024} introduces conservatism and relaxing this conservatism is an important open problem, next to being applicable to broader classes of nonlinear control systems, beyond (linear) reset controllers. 

To fill this gap, in this paper, we propose a novel framework based on the notion of dissipativity \cite{van_der_schaft_l2-gain_2017},\cite{Willems07} for computing scaled graphs {of possibly multivariable systems,} that reduces the aforementioned conservatism, and even leads to exact representations of scaled graphs in certain cases. To formulate the construction of (over-approximations of) scaled graphs, we exploit connections between LMIs, integral quadratic constraints (IQCs) \cite{megretski_system_1997}, and scaled graphs, which leads to  a convex optimization procedure in terms of solving sets of LMIs. Different from the approach in \cite{van_den_eijnden_scaled_2024}, we take the intersection of regions in the complex plane corresponding to input-output properties of the systems under study. We show that for specific (normal) LTI systems, our method leads to an exact representation of the scaled (relative) graph. This encouraging, non-conservative result motivates extension of the ideas and methods towards classes of nonlinear systems, including impulsive and piecewise linear systems. The latter two system classes are specific examples of hybrid systems \cite{goebel2012hybrid}, and arise in various relevant fields and engineering applications such as high-performance motion control systems \cite{van_den_eijnden_hybrid_2020}, \cite{Shi}, \cite{zaccarian_first_2005}, \cite{ZHANG2024106063}, non-smooth mechanical systems with impacts \cite{Leine2008}, and neural networks \cite{Samanipour}, to name but a few. 

We stress that the graphical nature of scaled graphs offers unique properties in terms of interpretability for system analysis and controller (re)design compared to, e.g., direct use of LMIs for analysis and controller design. For instance, robustness margins in feedback interconnections can be directly inferred from scaled graphs of the individual feedback components, and stabilizing controller classes can be visualized through scaled graph properties (in terms of their shape and location in the complex plane). These features further reinforce the importance of developing effective computational methods for scaled graphs, as is also acknowledged in \cite{chaffey_graphical_2023}, and to which we further contribute in this paper.

In line with the above, the main contributions of this paper are threefold: i) We establish formal connections between LMIs, IQCs and scaled graphs, allowing for an extension of the classical KYP-lemma for LTI systems \cite{rantzer_kalmanyakubovichpopov_1996}, ii) we use this novel result to formulate an effective LMI-based procedure for computing the scaled graph of LTI systems, and show that in certain cases this procedure is exact, and iii) we extend the framework towards impulsive and piecewise linear systems, and demonstrate the effectiveness through several numerical examples. In particular, we show how our new method significantly reduces conservatism and improves scaled graph computations for reset control systems as compared to the approach in \cite{van_den_eijnden_scaled_2024}. Initially, we focus on so-called \emph{soft} scaled graphs \cite{chen2025softhardscaledrelative}, which act on the $\mathcal{L}_2$-space. However, we will also show that our framework allows for computing \emph{hard} scaled graphs, acting on the extended $\mathcal{L}_2$-space, see \cite{chen2025softhardscaledrelative}. Soft and hard scaled graphs are closely related to soft and hard IQCs \cite[Section IV]{megretski_system_1997}, as we will highlight in this paper. 

The remainder of this paper is organized as follows. In Section~\ref{sec:prelim} we provide preliminaries on signals and systems, scaled graphs, and IQCs. Section~\ref{sec:mainidea} presents the main idea underlying our scaled graph computation method. Our first main result is presented in Section~\ref{sec:LTI} in the form of an explicit optimization procedure for LTI systems, along with an exactness result. The extension towards nonlinear system classes forms the second main result in this paper and is presented in Section~\ref{sec:NL}. We will comment on computations of hard scaled graphs in Section~\ref{sec:hardSG}. The main conclusions are provided in Section~\ref{sec:concl}.

\section{Preliminaries}
\label{sec:prelim}
Let $\mathbb{R}$, $\mathbb{R}_{\geq 0}$, and $\mathbb{R}_{> 0}$ denote the field of real numbers, non-negative real numbers, and strictly positive real numbers, and $\mathbb{C}$ the set of complex numbers. We denote the real and imaginary parts of a complex number $z\in \mathbb{C}$ by $\textup{Re}\left\{z\right\}$ and $\textup{Im}\left\{z\right\}$, respectively. For a nonzero complex number $z \in \mathbb{C}$ in polar form $|z|e^{j \angle z}$, the magnitude is denoted by $|z|$ and the angle is denoted by $\angle z$. The complex conjugate is denoted by $\bar{z}$, i.e., $\bar{z} = |z|e^{-j\angle z}$. Let $\mathbb{C}_{\geq 0}$ and $\mathbb{C}_{>0}$ denote the closed and open right-half complex plane, respectively, i.e., $\mathbb{C}_{\geq0} := \{z\in\mathbb{C} \mid \textup{Re}\{z\}\geq0\}$ and $\mathbb{C}_{>0} := \{z\in\mathbb{C} \mid \textup{Re}\{z\}>0\}$. Moreover, let $\mathbb{C}_+ := \{z\in\mathbb{C}\mid\textup{Im}\{z\}\geq0\}$, $\mathbb{C}_- := \{z\in\mathbb{C}\mid\textup{Im}\{z\}\leq0\}$, and let $\mathbb{D} := \{ z \in \mathbb{C} \mid |z| \leq 1 \}$ denote the closed unit disk. The closure of a set ${S}\subset \mathbb{C}$ is denoted by $\overline{{S}}$. The distance between two sets ${S}_1,S_2\in \mathbb{C}$ is denoted $\textup{dist}({S}_1,{S}_2)=\inf_{a\in {S}_1,\;b\in{S}_2}|a-b|$. A matrix $A \in \mathbb{R}^{n\times n}$ is said to be Hurwitz, if all its eigenvalues have strictly negative real part. The sets of $n$-by-$n$ real symmetric matrices are denoted by $\mathbb{S}^{n}=\{P\in\mathbb{R}^{n\times n}\mid P=P^\top\}$. For $P \in \mathbb{S}^n$, we use $P\succ 0$ (resp. $P\succeq 0$) to indicate that $P$ is positive definite (resp. positive semi-definite), i.e., $x^\top P x>0$ for all $x\in\mathbb{R}^{n}\setminus\{0\}$ (resp. $x^\top P x \geq 0$ for all $x\in \mathbb{R}^n$). Similar conventions are used for negative (semi-)definite matrices. The set of $n$-by-$n$ real symmetric matrices with non-negative elements is denoted by $\mathbb{S}_{\geq 0}^n$. Inequality symbols $>,\geq, <, \leq$ applied to vectors or matrices are understood component wise. The Kronecker product is denoted by $\otimes$ and $\textup{diag}(a_1,a_2,...,a_n)$ denotes the $n$-by-$n$ matrix with elements $a_1,a_2,...,a_n$ on its diagonal and $0$ elsewhere.

\subsection{Signals and systems}
We denote the space of vector-valued square-integrable functions with non-negative time-support by 
$$\mathcal{L}_2^n = \left\{u: \mathbb{R}_{\geq 0}\to\mathbb{R}^n \mid \int_{0}^\infty u(t)^\top u(t) dt < \infty\right\},$$ 
where the superscript $n$ is dropped when the dimension is clear from the context. For signals $u,y\in \mathcal{L}_2^n$, the inner product and induced norm are {given by}
$$\langle u,y\rangle = \int_{0}^\infty u(t)^\top y(t) dt \:\textup{ and }\:\|u\| = \sqrt{\langle u,u\rangle}.$$ 
Associated with the normed $\mathcal{L}_2^n$-space is the extended $\mathcal{L}_2$-space, given by
\begin{equation*}
      \mathcal{L}_{2e}^n = \left\{u:\mathbb{R}_{\geq 0} \to \mathbb{R}^n \mid P_Tu \in \mathcal{L}_2 \ \ \textup{for all } T \geq 0\right\},
\end{equation*}
where for $T\geq 0$, $P_T$ is the truncation operator on a signal $u:\mathbb{R}_{\geq 0}\to \mathbb{R}^n$ such that $P_Tu :\mathbb{R}_{\geq 0}\to \mathbb{R}^n$ with $(P_Tu)(t)=u(t)$ when $t\in [0,T]$, and $(P_Tu)(t)=0$ when $t> T$.

In this paper, $H$ denotes a system that maps inputs $u \in \mathcal{L}_{2}$ to outputs $y \in \mathcal{L}_{2}$. As a shorthand notation, we write $y\in H(u)$ to indicate all possible output trajectories related to the input $u$ when applied to the system $H$. When having unique solutions we write $y=H(u)$. We assume that the system maps zero input signals to zero output signals, i.e., $H(0)=0$. We will refer to systems mapping $\mathcal{L}_2$ to $\mathcal{L}_2$ as stable systems (in $\mathcal{L}_2$-sense). We call these systems bounded, if, in addition, $\sup_{u \in \mathcal{L}_2\setminus \left\{0\right\}}\sup_{y\in H(u)}\|y\|/\|u\|<\infty$. A system is said to be causal if $P_Ty \in P_TH(P_Tu)$ for all $T\geq 0$ and all $u \in \mathcal{L}_2$. We will only consider \emph{square} systems, i.e., systems with an equal number of inputs and outputs. We assume that $H$ can be modeled in a state-space framework with associated states $x \in \mathbb{R}^m$ representing the internal (latent) variables. To keep generality at this point, we do not commit to a specific state-space structure yet, but we will specify structure when needed.

\subsection{Scaled graphs}\label{sec:sg}
Let a system $H:\mathcal{L}_2\rightrightarrows \mathcal{L}_2$ be given. For signals $u \in \mathcal{U}\subset \mathcal{L}_2$ and $y \in H(u)$ define the \emph{gain} $\rho(u,y) \in [0, \infty]$ as
    \begin{align}\label{eq:gain}
    \rho(u,y) := \frac{\|y\|}{\|u\|},
    \end{align}
{if $u\neq 0$, and $\rho(u,y)$ is undefined if $u=0$. The \emph{phase} $\theta(u,y) \in [0,\pi]$ is defined as} 
\begin{align}\label{eq:phase}
    \theta(u,y):=&\arccos \frac{\langle u,y\rangle}{\|u\|\|y\|},
\end{align}
if $u,y\neq 0$, $\theta(u,y):=0$ if $u\neq 0, y=0$, and $\theta(u,y)$ is undefined if $u=0$. Note that by assumption $u=0$ yields $y=H(0)=0$, which in, e.g., feedback stability analysis leads to a trivial situation.   
The scaled graph of $H$ over inputs $u\in \mathcal{U}$ is defined as \cite{chaffey_graphical_2023}, \cite{chen2025softhardscaledrelative}
\begin{equation}\label{eq:SG}
    \textup{SG}_\mathcal{U}(H) := \left\{\rho(u,y)e^{\pm j\theta(u,y)} \mid u \in \mathcal{U}\setminus \left\{0\right\}, \:y \in H(u)\right\}.
\end{equation}
When $\mathcal{U}=\mathcal{L}_2$ we write $\textup{SG}(H)$. The inverse of $\textup{SG}(H)$ is obtained by swapping the role of the input
and output of $H$. In particular, we define the gain of the inverse scaled graph as $\rho(y,u)=1/\rho(u,y)$ for all $u,y\neq 0$ and $\rho(y,u):=\infty$ for $u\neq 0, y=0$, and the phase as $\theta(y,u) = \theta(u,y)$ for all $u \neq 0$. For $u=0$ the gain and phase are undefined. The inverse scaled graph of $H$ over signals $u\in \mathcal{U}$ is then given by 
\begin{equation}\label{eq:inv_SG}
\textup{SG}_\mathcal{U}^\dagger(H):=\left\{\rho(y,u)e^{\pm j\theta(y,u)}\mid u \in \mathcal{U}\setminus\left\{0\right\}, y\in H(u)\right\}.
\end{equation}
As before, in case $\mathcal{U}=\mathcal{L}_2$ we write $\textup{SG}^\dagger(H)$. The inverse scaled graph $\textup{SG}^\dagger(H)$ plays an important role in studying stability of the feedback interconnection depicted in Fig.~\ref{fig:FB}, as stated in the next theorem.
\tikzset{
    block/.style = {draw, rectangle,
        minimum height=1cm,
        minimum width=1cm},
    input/.style = {coordinate,node distance=1cm},
    output/.style = {coordinate,node distance=1cm},
    arrow/.style={draw, -latex,node distance=2cm},
    sum/.style = {draw, circle, node distance=.5cm, inner sep=.3em},
}      
\begin{figure}[htb]
\centering
\begin{tikzpicture}[auto, node distance=1cm,>=latex']
\centering
\node [block] (H1) {$H_1$};
\node [sum, left=0.7 of H1] (s1) {};         \node [block, below=.6 of H1] (H2) {$H_2$};
\node [sum, right=.7 of H2] (s2) {};

\draw[->] (s1) --node [xshift = -2, yshift = 2 , pos=0.5] {$u_1$}  (H1);
\draw[->] (s2) --node [xshift = 2, yshift = -2 , pos=0.5] {$u_2$}  (H2);

\draw[->] (H1) -| node [pos=0.92] {\small $+$} (s2);
\draw[->] (H2) -| node [pos=0.92] {\small $-$} (s1);

\node [input, left=.8 of s1] (dummy1){};
\draw[->] (dummy1) --node [xshift = -2, yshift = 2 , pos=0.5] {$w_1$}(s1);

\node [input, right=0.8 of s2] (dummy2){};
\draw[->] (dummy2) -- node [xshift = 2, yshift = -2 , pos=0.5] {$w_2$}(s2);

\node [input, left=1.85 of H2] (dummy3){};
\draw[->] (H2) --node [xshift = 2, yshift = -2 , pos=0.8] {$y_2$}(dummy3);

\node [input, right=1.85 of H1] (dummy4){};
\draw[->] (H1) --node [xshift = -2, yshift = 2 , pos=0.8] {$y_1$}(dummy4);

\end{tikzpicture}
 \caption{Feedback interconnection.}
 \label{fig:FB}
\end{figure}
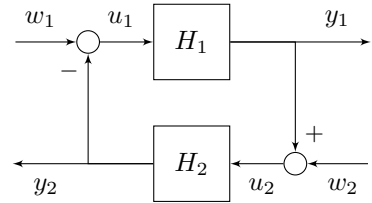
\begin{thm}\textnormal{\cite[Theorem 2]{chen2025softhardscaledrelative}}\label{th:SGfb}
Consider a pair of causal, stable systems $H_1:\mathcal{L}_2\to\mathcal{L}_2$ and $H_2:\mathcal{L}_2\to\mathcal{L}_2$ placed in negative feedback interconnection as depicted in Fig.~\ref{fig:FB}. Suppose that the interconnection of $H_1$ and $\tau H_2$ is well-posed\footnote{A feedback interconnection is well-posed if, given input signals in $\mathcal{L}_2$, there exist output signals in $\mathcal{L}_{2e}$ depending causally on the inputs, see \cite{freeman_role_2022}.} for all $\tau \in (0,1]$. If there exists $r >0$ such that for all $\tau \in (0,1]$ 
\begin{equation}\label{eq:cond}
   \textup{dist}(\textup{SG}^\dagger(-H_1), \textup{SG}(\tau H_2))\geq  r,
\end{equation}
then inputs $w = [w_1^\top, w_2^\top]^\top \in \mathcal{L}_2$ are mapped to outputs $y = [y_1^\top, y_2^\top]^\top\in \mathcal{L}_2$, and $\|y\|/\|w\|\leq 1/r$. \hfill $\ulcorner$
\end{thm}

The attractive aspect of Theorem~\ref{th:SGfb} is the fact that condition \eqref{eq:cond}, and thus feedback stability can be checked graphically, and is reminiscent of a classical Nyquist stability test for LTI systems. Moreover, the distance between the scaled graphs directly provides an indication of performance (in terms of the $\mathcal{L}_2$-gain), as well as robustness measures (in terms of distances of sets). These properties underscore the strengths of the scaled graph framework, and highlight the need for accurately estimating scaled graphs of systems.\footnote{Note that the result in Theorem~\ref{th:SGfb} also holds when using over-approximations of $\textup{SG}^\dagger(-H_1)$ and $\textup{SG}(H_2)$, in which case \eqref{eq:cond} is clearly satisfied, but at the cost of conservatism.} 

The definition of the scaled graph in \eqref{eq:SG} is also referred to as a \emph{soft} scaled graph and Theorem~\ref{th:SGfb} is referred to as soft scaled graph separation \cite{chen2025softhardscaledrelative}. Soft scaled graphs deal with trajectories in $\mathcal{L}_2$. In contrast, \emph{hard} scaled graphs introduced in \cite{chen2025softhardscaledrelative} deal with trajectories in $\mathcal{L}_{2e}$, and can accommodate input-output behaviour of possibly unbounded systems such as LTI integrators, see \cite{chen2025softhardscaledrelative}. In the following sections, we primarily focus on computing soft scaled graphs; definitions and computations for hard scaled graphs will be postponed to Section~\ref{sec:hardSG}.

\subsection{Hyperbolic Geometry} \label{sec:hyperbolic_geometry}
Since we will make use of elements of hyperbolic geometry \cite{chaffey_graphical_2023}, \cite{huang_scaled_2020}, \cite{pates_scaled_2021}, we recall here some of the necessary concepts. Two models of hyperbolic geometry will play a key role: the Poincaré upper half-plane model and the Beltrami-Klein disk model of hyperbolic space. First, we recall the relevant notions associated with the Poincaré upper half-plane model. In this model, the geodesic $\text{ARC}_{\text{min}}(z_1,z_2)$, between two points $z_1,z_2\in\mathbb{C}_{+}$ is the section of the circle centered on the real axis through $z_1,z_2$, with endpoints $z_1$ and $z_2$, that lies in $\mathbb{C}_+$. Under the Poincaré upper half-plane model, we adopt the definitions of hyperbolic convexity (or h-convexity) and the h-convex hull from \cite[Definition 6]{chaffey_graphical_2023}.
\begin{defn}\textnormal{\cite[Definition 6]{chaffey_graphical_2023}} A set $S\subset\mathbb{C}_{+}$ is said to be h-convex, if
	\begin{align} 
		z_1,z_2\in S \Rightarrow  \textup{ARC}_{\textup{min}}(z_1,z_2)\in S.
		\end{align} 
        Given a set of points $P \subset \mathbb{C}_+$, the h-convex hull of $P$, denoted by h-hull($P$), is the smallest h-convex set containing $P$.
\end{defn}
We denote the standard Euclidean convex hull of a set of points $P \subset \mathbb{C}$ as hull($P$). 


Next, we recall the Beltrami-Klein mapping $f_{BK}:\mathbb{C}\rightarrow\mathbb{D}$, which is given by
\begin{align}
	f_{BK}(z)=\frac{(\bar{z}-j)(z-j)}{1+\bar{z}z}=\frac{(\bar{z}-j)(z-j)}{1+|z|^2}, \label{eq:bk_map}
\end{align} 
 which maps geodesics $\textup{ARC}_{\textup{min}}(z_1,z_2)$ to straight line segments. As discussed in \cite{pates_scaled_2021}, the Beltrami-Klein map in \eqref{eq:bk_map} is bijective on $\mathbb{C}_+$ as a maping from $\mathbb{C}_+\to \mathbb{D}$, but not on $\mathbb{C}$, as for all $z\in\mathbb{C}$, $f_{BK}(z)=f_{BK}(\bar{z})$. This motivates the inverse map $f_{BK}^{-1}:\mathbb{D}\rightrightarrows\mathbb{C}$, given by
\begin{align} 
	f_{BK}^{-1}(z) = \left\{\frac{\textup{Im}\{z\}\pm j\sqrt{1-|z|^2}}{\textup{Re}\{z\}-1}\right\},
\end{align} 
such that $f_{BK}^{-1}(f_{BK}(z))=\{z,\bar{z}\}$ for all $z\in\mathbb{C}$, see \cite{pates_scaled_2021}. Note that $f_{BK}^{-1}$ is set-valued injective on $\mathbb{D}$, in the sense that for $a,b \in \mathbb{D}$, $f^{-1}_{BK}(a) \cap f^{-1}_{BK}(b) \neq \emptyset$ implies $a = b$.

The function $f_{BK}$ maps arcs centered on the real axis to Euclidean straight line segments in the unit disk, since in the Beltrami-Klein model geodesics are Euclidean straight lines. Consequently, h-convexity in the Beltrami-Klein model coincides with Euclidean convexity \cite{pates_scaled_2021}.
This connection between hyperbolic convexity and Euclidean convexity results, for a given ${S}\subset\mathbb{C}_+$, in
\begin{align} 
\textup{h-hull}({S})=f_{BK}^{-1}(\textup{hull}(f_{BK}({S}))).
\end{align} 


\subsection{Dissipativity and IQCs}\label{sec:dissp}
The main aim of this paper is the development of methods for computing (a tight over-approximation of) the scaled graph in \eqref{eq:SG} of a system $H$. The machinery that we will exploit for this purpose is based on the well-known notion of dissipativity, see \cite{van_der_schaft_l2-gain_2017,Willems07}. Recall that we assume that $H$ can be modeled within a state-space framework, where states $x: \mathbb{R}_{\geq 0}\to \mathbb{R}^m$ represent latent (internal) variables in the mapping from $u$ to $y$.

\begin{defn}\label{def:dissipativity}
A system $H:\mathcal{L}_{2}^n \rightrightarrows \mathcal{L}_{2}^n$ with corresponding state trajectories $x:\mathbb{R}_{\geq 0}\to \mathbb{R}^m$ is said to be dissipative with respect to the supply rate ${s}:\mathbb{R}^n\times \mathbb{R}^n \to \mathbb{R}$, if there exists a function $W:\mathbb{R}^m \rightarrow\mathbb{R}$ called the storage function that satisfies $W(0)=0$ and \begin{align}\label{eq:def_dissipativity}
W(x(T))\leq W(x(0))+\int_{0}^{T} {s}(u(t),y(t))\;dt
\end{align} 
for all $T\geq 0$, and all possible trajectories $(u,x,y)$ generated by $H$.
\end{defn}
Note that in the above definition of dissipativity we do not require $W$ to be non-negative. If we assume that $\lim_{t\to \infty} x(t) = 0$ (which is the case if, for example both $x\in \mathcal{L}_2$ and $\dot{x} \in \mathcal{L}_2$ \cite[p. 237]{Desoer}), and the system is initially at rest so that $x(0) = 0$, then, using $W(0) = 0$, the inequality in \eqref{eq:def_dissipativity} implies
\begin{equation}\label{eq:diss_inf}
 0\leq \int_{0}^\infty s(u(t),y(t))dt.
\end{equation}
When restricting $W(x)\geq 0$ for all $x \in \mathbb{R}^m$, the inequality holds without the extra assumption on $x$. In fact, in that case we have $0\leq \int_0^T  s(u(t),y(t))dt$ for all $T\geq 0$. 

In the remainder of this paper, we will mainly work with quadratic supply rates of the form
\begin{align}\label{eq:canonical_supply_rate}
s(u,y) =\begin{bmatrix} y \\ u \end{bmatrix}^\top \left(\Pi \otimes I_n \right)\begin{bmatrix} y \\ u \end{bmatrix},  \textup{ with } \Pi = \begin{bmatrix} a & b\\b & c \end{bmatrix},
\end{align}
and where $a,b,c \in \mathbb{R}$. In this case, \eqref{eq:diss_inf} reduces to an IQC \cite{megretski_system_1997}. The link between dissipativity and IQCs is well established, see, e.g., \cite{scherer_linear_2015,van_der_schaft_l2-gain_2017}. This paper further details the link between dissipativity, IQCs, and scaled graphs.

\section{Over-approximation of SGs}
\label{sec:mainidea}

In this section, the main idea and approach for approximating the scaled graph of a system $H$ is presented. To this end, we first establish connections between IQCs and scaled graphs. 

\subsection{ A connection between IQCs and scaled graphs}
Consider the IQC
\begin{align} \label{eq:Supply_inequality}
\int_{0}^{\infty}\begin{bmatrix} y(t) \\ u(t) \end{bmatrix}^\top \left(\Pi \otimes I_n\right)\begin{bmatrix} y(t) \\ u(t) \end{bmatrix} dt \geq 0.
\end{align} 
In \cite[Lemma 1]{van_den_eijnden_scaled_2024} a link between IQCs of the form \eqref{eq:Supply_inequality} and the SG in \eqref{eq:SG} has been established for the single-input single-output case, i.e., $n=1$ and thus, $I_n=1$. Below we slightly extend this result to accommodate the multivariable case.

\begin{lem}\label{lem:connection_to_SG}
A stable system $H:\mathcal{L}_{2} \rightrightarrows \mathcal{L}_{2}$ satisfies the IQC in \eqref{eq:Supply_inequality} for all $u\in \mathcal{U}$ and $y \in H(u)$ if and only if $\textup{SG}_{\mathcal{U}}(H) \subset \mathcal{S}(\Pi)$, with
\begin{align} 
	\mathcal{S}(\Pi) = \left\{ z \in\mathbb{C} \left|  \begin{bmatrix}
		z \\ 1
	\end{bmatrix}^* \Pi \begin{bmatrix}
		z \\ 1
	\end{bmatrix} \geq 0 \right.  \right\}.\label{eq:S_PI_canonical}
\end{align} 
When $\textup{det}(\Pi)<0$, then for:
\begin{enumerate}
    \item $a=0$, $b >0$ the region $\mathcal{S}(\Pi)$ describes the half-plane $\textup{Re}(z) \geq -c/2b$;
    \item $a=0$, $b <0$ the region $\mathcal{S}(\Pi)$ describes the half-plane $\textup{Re}(z)\leq-c/2b$;
    \item $a>0$, the region $\mathcal{S}(\Pi)$ describes the exterior of the disc with center $-b/a$ and radius $\sqrt{b^2/a^2 - c/a}$;
    \item $a<0$, the region $\mathcal{S}(\Pi)$ describes the interior of the disc with center $-b/a$ and radius $\sqrt{b^2/a^2 - c/a}$.
\end{enumerate}
\end{lem}
\begin{pf}
Note that the IQC in \eqref{eq:Supply_inequality} can be written as
\begin{equation}\label{eq:ineq1}
    0 \leq a\|y\|^2 + 2b\langle u,y\rangle + c\|u\|^2. 
\end{equation}
Denote $z(u,y) = \rho(u,y)e^{\pm j\theta(u,y)}$ with $\rho(u,y)$ and $\theta(u,y)$ the gain and phase defined in \eqref{eq:gain} and \eqref{eq:phase}, respectively. As such, we find $|z(u,y)| = \rho(u,y)=\|y\|/\|u\|$ and $\textup{Re}\left\{z(u,y)\right\} = \rho(u,y)\cos(\theta(u,y)) = \langle u,y\rangle/\|u\|^2$. Then, dividing \eqref{eq:ineq1} by $\|u\|^2$, we find
\begin{equation}\label{eq:ineq2}
\begin{split}
   &a|z(u,y)|^2+2b\textup{Re}\left\{z(u,y)\right\}+c\\
   &\qquad=\begin{bmatrix}
       z(u,y)\\1
   \end{bmatrix}^*\begin{bmatrix}
       a & b\\b &c
   \end{bmatrix}\begin{bmatrix}
       z(u,y)\\1
   \end{bmatrix}\geq 0,
   \end{split}
\end{equation}
which shows that $\textup{SG}_{\mathcal{U}}(H) \subset \mathcal{S}(\Pi)$. 

If $a = 0$, case (1) and (2) follow immediately from \eqref{eq:ineq2}. To show case (3) and case (4), we let $x = \textup{Re}\left\{z\right\}$ and $y = \textup{Im}\left\{z\right\}$, and rewrite the inequality in \eqref{eq:ineq2} as
\begin{subequations}
\begin{align}
    \left(x+\frac{b}{a}\right)^2 + y^2 + \frac{c}{a} - \frac{b^2}{a^2} & \geq 0\quad \textup{if } a >0, \label{eq:disc1}\\
    \left(x+\frac{b}{a}\right)^2 + y^2 + \frac{c}{a} - \frac{b^2}{a^2} & \leq 0 \quad \textup{if } a <0.\label{eq:disc2}
    \end{align}
\end{subequations}
Suppose $\textup{det}(\Pi) = ac-b^2 <0$, such that we have $c/a-b^2/a^2 <0$. Then, \eqref{eq:disc1} and \eqref{eq:disc2} specify the exterior respectively interior of a disc centered at $-b/a$ and having radius $\sqrt{b^2/a^2-c/a}$. \hfill $\blacksquare$
\end{pf}
Lemma~\ref{lem:connection_to_SG} provides a direct connection between an IQC for $H$ and an over-approximation of $\textup{SG}_{\mathcal{U}}(H)$. This suggests that the more IQCs of the form \eqref{eq:Supply_inequality} are available for the system $H$, the more over-approximating regions we can find. Taking the intersection of these regions then allows to tighten the over-approximation of $\textup{SG}_\mathcal{U}(H)$. This observation is at the core of our SG approximation method, which we discuss in detail next. 

\subsection{Approximation method}\label{sec:approximation_1}
Given a system $H$ taking inputs $u \in \mathcal{U}$ and producing outputs $y \in \mathcal{L}_2$, we define the set of matrices 
\begin{multline}\label{eq:PI_set}
\mathbf{\Pi}(H) = \left\{ \Pi \in \mathbb{S}^2 \;\middle|\; \det(\Pi) < 0 \text{ and } \right. \\
\left. \qquad \eqref{eq:Supply_inequality} \text{ holds for all } u \in \mathcal{U},\, y \in H(u) \right\}.
\end{multline}
The following result provides an over-approximation of $\textup{SG}_\mathcal{U}(H)$ exploiting multiple (IQC generated) regions.
\begin{thm}\label{th:general_SG_approx}
    Consider a system $H:\mathcal{L}_2\rightrightarrows\mathcal{L}_2$ and suppose $\mathbf{\Pi}(H) \neq \emptyset$. Let $\tilde{\mathbf{\Pi}} \subset \mathbf{\Pi}(H)$. Then
        \begin{align} 
\textup{SG}_\mathcal{U}(H)\subseteq \bigcap_{\Pi\in\tilde{\mathbf{{\Pi}}}}\mathcal{S}(\Pi).
\end{align} 
\end{thm}
\begin{pf}
For any $\Pi\in \tilde{\mathbf{\Pi}}$ it holds that $\Pi \in \mathbf{\Pi}(H)$ and thus it follows that $\textup{SG}_\mathcal{U}(H) \subset \mathcal{S}(\Pi)$. Since $\Pi \in \tilde{\mathbf{\Pi}}$ is arbitrary, we have $\textup{SG}_\mathcal{U}(H)\subseteq \bigcap_{\Pi\in\tilde{\mathbf{{\Pi}}}}\mathcal{S}(\Pi)$. \hfill $\blacksquare$
\end{pf}

\begin{figure}[htb]
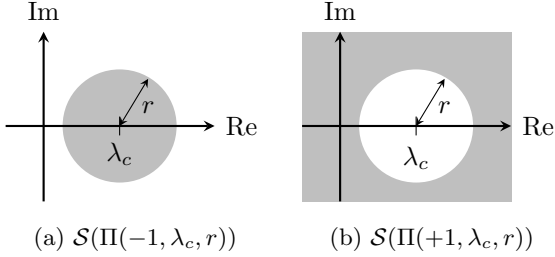

	\centering
	\begin{subfigure}{0.45\linewidth}
		{\ctikzfig{Figures/General_complex_circle}}
		\caption{$\mathcal{S}(\Pi(-1,\lambda_c,r))$}
        \label{fig1a}
	\end{subfigure}
	\begin{subfigure}{0.45\linewidth}
		{\ctikzfig{Figures/Exterior_General_complex_circle}}
		\caption{$\mathcal{S}(\Pi(+1,\lambda_c,r))$}
        \label{fig1b}
	\end{subfigure}
	
	\caption{Regions defined by $\mathcal{S}(\Pi(\sigma,\lambda_c,r))$ indicated in grey for (a) $\sigma = -1$ and (b) $\sigma=+1$.}
	\label{fig:in-out_circle}
\end{figure}

Note that in Theorem~\ref{th:general_SG_approx} we work with $\textup{det}(\Pi)<0$. When $\textup{det}(\Pi) \geq 0$, the region $\mathcal{S}(\Pi)$ in \eqref{eq:S_PI_canonical} yields an empty set or the entire complex plane, see also the discussion in \cite[Section IV]{Iwasaki}. Therefore, in principle, the interesting case to consider is $\textup{det}(\Pi)<0$ leading to (interior/exterior) circular regions and half-planes as shown in Lemma~\ref{lem:connection_to_SG}. In this context, a useful parametrization for the matrix $\Pi$ to consider is given by
\begin{equation}\label{eq:PI_parametrization}
    \Pi(\sigma, \lambda, r) = \sigma\begin{bmatrix}
        1 & -\lambda_c \\
        -\lambda_c & \lambda_c^2-r^2
    \end{bmatrix},
\end{equation}
with $\sigma \in \left\{-1, +1\right\}$, $\lambda_c \in \mathbb{R}$, and $r \in \mathbb{R}_{> 0}$. Note that $\textup{det}(\Pi(\sigma,\lambda_c,r)) = -r^2 <0$. As such, $\Pi(-1,\lambda_c,r)$ represents the interior of a disk in the complex plane centered on the real axis at $\lambda_c$ and with radius $r$, and $\Pi(+1,\lambda_c,r)$ represents the exterior of such a disk. Examples of the generated regions are provided in Fig.~\ref{fig:in-out_circle}.

\begin{figure}[t]
		\centering
		\includegraphics[width=0.7\linewidth]{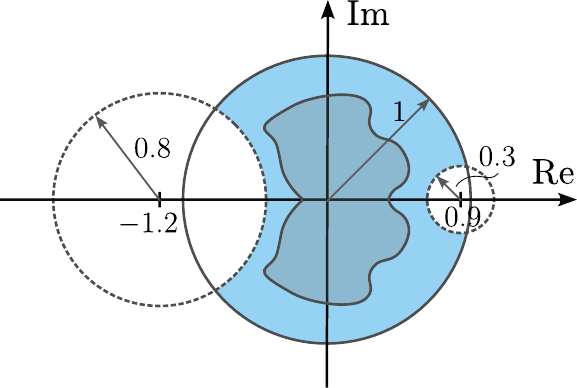}
		\caption{Illustration of $\textup{SG}(H)$ (dark blue region) and over-approximation (light blue region), generated by the matrices $\Pi(-1,0,1)$ (interior of solid circle), $\Pi(+1,-1.2,0.8)$ (exterior of left dashed circle), and $\Pi(+1,0.9,0.3)$ (exterior of right dashed circle).}
        \label{fig:simple_circle_example}
	\end{figure}

To generate an over-approximation of the scaled graph of a system $H$ we search for different parameter combinations $(\sigma, \lambda_c, r)$ such that matrices $\Pi$ of the form \eqref{eq:PI_parametrization} satisfy \eqref{eq:Supply_inequality}. Collecting all resulting matrices in a set $\tilde{\mathbf{\Pi}}\subseteq \mathbf{\Pi}(H)$ leads via Theorem~\ref{th:general_SG_approx} to an over-approximation of $\textup{SG}_{\mathcal{U}}(H)$ generated by all matrices in $\tilde{\mathbf{\Pi}}$. A sketch of such an over-approximation is shown in Fig.~\ref{fig:simple_circle_example} in blue. Note that using the exterior of a disk allows for obtaining non-convex regions (in Euclidean sense). As we will see later, this offers a significant benefit in constructing tight SG over-approximations as compared to, e.g., the method in \cite{van_den_eijnden_scaled_2024}, which only exploits interior regions.

The parameterization of the matrix $\Pi$ in \eqref{eq:PI_parametrization} also allows for computations of over-approximations of the inverse scaled graph $\textup{SG}^\dagger(H)$ as defined in \eqref{eq:inv_SG}.
\begin{lem}
     Consider a system $H:\mathcal{L}_2\rightrightarrows\mathcal{L}_2$ and suppose that $\textup{SG}(H) \subset \mathcal{S}(\Pi(\sigma,\lambda_c,r))$, with $\sigma\in\{-1,+1\}$, $\lambda_c\in\mathbb{R}$ and $r\in\mathbb{R}_{>0}$. Let $\lambda'_c = \lambda_c/(\lambda_c^2-r^2)$ and $r' = r/|\lambda_c^2-r^2|$.
\begin{enumerate}
    \item If $0 \not\in\mathcal{S}(\Pi(\sigma,\lambda_c,r))$, then $\textup{SG}^\dagger(H) \subset \\\mathcal{S}(\Pi(-1,\lambda_c{'} ,r{'}))$;
        \item If $0 \in \mathcal{S}(\Pi(\sigma,\lambda_c,r))$ and $|\lambda_c| \neq r$, then $\textup{SG}^\dagger(H)\subset \mathcal{S}(\Pi(+1,\lambda_c{'} ,r{'}))$;
        \item If $0\in \mathcal{S}(\Pi(\sigma,\lambda_c,r))$, $|\lambda_c| = r$, and $\lambda_c>0$, then $\textup{SG}^\dagger(H)\subset\left\{z\in\mathbb{C}\mid -\sigma\textup{Re}\{z\}\geq \frac{-\sigma}{\lambda_c+r}\right\}$;
    \item If $0\in \mathcal{S}(\Pi(\sigma,\lambda_c,r))$, $|\lambda_c| = r$, and $\lambda_c<0$, then $\textup{SG}^\dagger(H)\subset\left\{z\in\mathbb{C}\mid-\sigma\textup{Re}\{z\}\leq \frac{\sigma}{\lambda_c-r}\right\}$.
\end{enumerate}
\end{lem}
\begin{pf}
The proof follows from standard geometric reasoning, see, e.g., \cite{needham1997visual}. A brief sketch is provided here. The region $\mathcal{S}(\Pi(\sigma, \lambda_c, r))$ denotes the interior ($\sigma = -1$) or exterior ($\sigma = +1$) of the circle $\mathcal{C}(\lambda_c, r)$, centered at $\lambda_c$ with radius $r$, recall Fig.~\ref{fig:in-out_circle}. 

Consider the inversion map \( z \mapsto ({1}/{|z|})e^{\pm j\angle z} \), also known as inversion in the unit circle. If \( 0 \notin \mathcal{C}(\lambda_c, r) \), i.e., \( |\lambda_c| \ne r \), this map transforms \( \mathcal{C}(\lambda_c, r) \) into another circle \( \mathcal{C}(\lambda_c', r') \). If \( 0 \notin \mathcal{S}(\Pi(\sigma, \lambda_c, r)) \), the inversion maps the exterior/interior of the original circle to the \emph{interior} of $\mathcal{C}(\lambda_c',r')$, yielding {case (1)}. If \( 0 \in \mathcal{S}(\Pi(\sigma, \lambda_c, r)) \), the origin maps to infinity, which lies in the \emph{exterior} of $\mathcal{C}(\lambda_c',r')$, leading to {case (2)}, see \cite[Chap.3]{needham1997visual} for details.

If \( 0 \in \mathcal{C}(\lambda_c, r) \), and $\lambda_c >0$, then \( \lambda_c = r \), and inversion maps the circle to a vertical line intersecting the real axis at \( 1/(\lambda_c + r) \). Points $z\in  \mathcal{S}(\Pi(-1, \lambda_c, r))$ map to the half-plane to the right of this line, and points  $z\in \mathcal{S}(\Pi(+1, \lambda_c, r))$ map to the left of this line, leading to {case (3)}. {Case (4)} follows by symmetry. \hfill $\blacksquare$
\end{pf}

\begin{rem}\label{rem:PIpass}
   Specific choices for $\Pi$ that may not be covered by \eqref{eq:PI_parametrization} (as $\lambda_c$ and $r$ take finite values) are given by
   \begin{equation}\label{eq:PIpasive}
       \Pi = \begin{bmatrix}
           0 & 1 \\
           1 & c
       \end{bmatrix}, \:\textup{ and } \:\Pi = \begin{bmatrix}
           0 & -1 \\
           -1 & c
       \end{bmatrix}, 
   \end{equation}
   where $c \in \mathbb{R}$, which correspond to (shifted) half-planes in $\mathbb{C}$, see item (1) and (2) in Lemma~\ref{lem:connection_to_SG}. These matrices are specifically useful to characterize (anti-)passive systems, and can be included in the set of matrices \eqref{eq:PI_set} without problems. However, the scaled graph of bounded systems will always belong to a bounded region in $\mathbb{C}$, as by definition $\|y\|/\|u\|\leq \gamma :=\sup_{u \in \mathcal{L}_2\setminus \left\{0\right\}}\sup_{y\in H(u)}\|y\|/\|u\|<\infty$, such that the effect of including \eqref{eq:PIpasive} can be approximated by taking $\lambda_c$ and $r$ large. Note that matrices of the form in \eqref{eq:PIpasive} may still have computational benefits. Such matrices may be particularly useful in computing hard scaled graphs that can deal with unbounded systems such as integrators. We will comment on this later in Section~\ref{sec:hardSG}. 
\end{rem}

    The rationale outlined above applies broadly to various system types, including static, dynamic, continuous-time, and discrete-time systems. When dissipativity with respect to multiple quadratic supply rates of the form \eqref{eq:canonical_supply_rate} can be verified for a system $H$, for example by constructing appropriate storage functions $W$ as in \eqref{eq:def_dissipativity}, an over-approximation of $\textup{SG}(H)$ follows directly from Theorem~\ref{th:general_SG_approx}. However, identifying matrices $\Pi(\sigma,\lambda_c,r)$ (or the corresponding parameters) that ensure \eqref{eq:Supply_inequality} holds for a given (nonlinear) system is generally nontrivial. In the next section, we present a procedure based on LMIs to address this problem for LTI systems. We then extend the approach to specific classes of nonlinear systems in the sections thereafter.

\section{Results for LTI Systems}
\label{sec:LTI}
\subsection{A KYP-like result}
Consider an LTI system $H_L:\mathcal{L}_2\rightarrow\mathcal{L}_2$, which can be described in state-space form by
\begin{equation}
    H_L:\;\left\{\begin{aligned}
    \dot{x}(t) &= Ax(t)+Bu(t),\qquad x(0) = 0,\\
    y(t) &= Cx(t)+Du(t),
    \end{aligned}\right. \label{eq:LTI_canonical}
\end{equation}
with state $x(t) \in \mathbb{R}^m$, input $u(t) \in \mathbb{R}^n$ and output $y(t) \in \mathbb{R}^n$ all at time $t \in \mathbb{R}_{\geq 0}$, and matrices $(A,B,C,D)$ of suitable dimensions, such that \eqref{eq:LTI_canonical} is minimal. The transfer function matrix from $u$ to $y$ is given by
\begin{equation}
    H_L(s) = C(sI-A)^{-1}B+D, \qquad s \in \mathbb{C}.
\end{equation}
We call an LTI system \emph{normal} if $\bar{H}_L(s) H_L(s) = H_L(s) \bar{H}_L(s)$, $s \in \mathbb{C}$. 

In the next result, we make a connection between the scaled (relative) graph of an LTI system and its corresponding transfer function matrix. We point out here that the scaled graph of an LTI system is equivalent to its scaled \emph{relative} graph \cite{chaffey_graphical_2023}.

 \begin{thm}\textnormal{\cite[Theorem 1 $ii$)]{pates_scaled_2021}}
 	\label{th:LTI_bk_hull}
 	Given a normal, stable LTI system $H_L:\mathcal{L}_2\rightarrow\mathcal{L}_2$ with transfer function matrix $H_L(s)$ and spectrum 
      \begin{equation}\label{eq:spec_}
        \sigma(H_L) = \cup_{\omega \in \mathbb{R}}\left\{ \lambda \in \mathbb{C} \mid \det(H_L(j\omega)-\lambda I)=0\right\}.
    \end{equation}
 Then 
 	\begin{align} 
 		{\textup{SG}}(H_L)=f_{BK}^{-1}(\textup{hull}(f_{BK}(\sigma(H_L)))).
	\end{align}
 \end{thm} 
 \begin{pf}
 The result follows directly from Theorem 1 ii) and Example~2 in \cite{pates_scaled_2021}.\hfill $\blacksquare$
 \end{pf}
In the single-input single-output (SISO) case, the spectrum in \eqref{eq:spec_} boils down to the Nyquist diagram of $H_L$, that is $\textup{Nyq}(H_L)=\left\{H_L(j\omega)\mid \omega \in \mathbb{R}\right\}$, so that for SISO LTI systems $${\textup{SG}}(H_L) = f_{BK}^{-1}(\textup{hull}(f_{BK}(\textup{Nyq}(H_L)))).$$
Hence, the scaled graph of a normal LTI system can be determined exactly from its corresponding characteristic loci and Nyquist contours, see also \cite[Theorem 4]{chaffey_graphical_2023}. 

In addition to the above result that exploits transfer functions, we also provide a connection between scaled graphs and a state-space representation by combining scaled graphs with the Kalman-Yakubovich-Popov (KYP) lemma \cite{rantzer_kalmanyakubovichpopov_1996}.

\begin{thm}
	\label{th:LTI_approx}
Consider an LTI system $H_L$ of the form (\ref{eq:LTI_canonical}) and assume that $A$ is Hurwitz. Let 
\begin{align} \label{eq:thetaPi}
\Theta(\Pi) = \begin{bmatrix}
		C & D \\ 0 & I
	\end{bmatrix}^\top(\Pi\otimes I_n)\begin{bmatrix}
		C & D \\ 0 & I
	\end{bmatrix}, \textup{ with } \Pi \in \mathbb{S}^2.
\end{align} 
Consider the following statements:
\begin{enumerate}
    \item There exists a matrix $P\in \mathbb{S}^m$ that satisfies the LMI
	\begin{align} \label{eq:KYP_LMI}
		 \begin{bmatrix}
			A &B\\I&0
		\end{bmatrix}^\top\begin{bmatrix}
		0 & P \\ P &0
		\end{bmatrix}\begin{bmatrix}
			A &B\\I&0
		\end{bmatrix} - \Theta(\Pi) &\preceq 0. 
    \end{align} 
    \item The IQC in \eqref{eq:Supply_inequality} holds for all $u \in \mathcal{L}_2$, $y \in H_L(u)$. 
    \item $\textup{SG}(H_L) \subset \mathcal{S}(\Pi)$, with $\mathcal{S}(\Pi)$ in \eqref{eq:S_PI_canonical}. 
\end{enumerate}
Then, $(1)\implies (2)$,  $(1)\implies (3)$, and $(2)\Leftrightarrow (3)$. Under controllability of $(A,B)$, all statements are equivalent.
\end{thm}
\begin{pf}
$(1)\implies (2)$ is a standard connection between IQCs and LMIs \cite{megretski_system_1997}. $(2)\Leftrightarrow (3)$ is a consequence of Lemma~\ref{lem:connection_to_SG}. $(1)\implies(3)$ follows accordingly. Under the additional controllability assumption, equivalence follows directly from the KYP-lemma \cite{rantzer_kalmanyakubovichpopov_1996}. {\hfill $\blacksquare$}
\end{pf}

Note that Theorem~\ref{th:LTI_approx} works with $\mathcal{U}=\mathcal{L}_2$. Accounting for $\mathcal{U}\neq\mathcal{L}_2$ can be done through an extension of the generalized KYP lemma \cite{Iwasaki} in the same spirit as \cite{van_den_eijnden_scaled_2024}. Both Theorem~\ref{th:LTI_bk_hull} and Theorem~\ref{th:LTI_approx} provide a means to compute (over-approximations of) scaled graphs for LTI systems. Theorem~\ref{th:LTI_bk_hull} exploits transfer function matrices, which is specific to LTI systems, whereas Theorem~\ref{th:LTI_approx} includes a state-space representation, the latter which is also relevant for nonlinear systems. Theorem~\ref{th:LTI_approx} forms the basis for developing a systematic numerical procedure to compute an over-approximation of $\textup{SG}(H_L)$ (and $\textup{SG}_\mathcal{U}(H_L)$) that can be extended towards nonlinear systems as well, as we will show later. A detailed description of the computation procedure follows next.

\subsection{Numerical procedure}\label{sec:proc}
Let $\Lambda_i$ and $\Lambda_e$ be subsets of $\mathbb{R}$, possibly containing a finite number of elements, and define the matrix
\begin{equation}\label{eq:LMI_F}
\begin{split}
    & F(\sigma, \lambda_c, r, P) := \\
    &\qquad \begin{bmatrix}
			A &B\\I&0
		\end{bmatrix}^\top\begin{bmatrix}
		0 & P \\ P &0
		\end{bmatrix}\begin{bmatrix}
			A &B\\I&0
		\end{bmatrix} - \Theta(\Pi(\sigma,\lambda_c,r)).
        \end{split}
\end{equation}
We formulate the following two optimization problems.

\textbf{L1:} For each $\lambda_c \in \Lambda_i$, 
\begin{equation}
\begin{aligned}
&\underset{r,P}{\text{minimize}}        & & r^2\\
&\text{subject to} & & F(-1,\lambda_c,r,P)\preceq 0, \\
& & & P \in \mathbb{S}^n, r> 0.
\end{aligned}
\end{equation}
\textbf{L2:} For each $\lambda_c \in \Lambda_e$, 
\begin{equation}
\begin{aligned}
&\underset{r,P}{\text{maximize}}        & & r^2\\
&\text{subject to} & & F(+1,\lambda_c,r,P)\preceq 0, \\
& & & P \in \mathbb{S}^n, r> 0.
\end{aligned}
\end{equation}

In \textbf{L1}, we minimize the interior region of the disk in the complex plane centered at $\lambda_c \in \Lambda_i$ (recall the shaded region in Fig.~\ref{fig1a}), whereas in \textbf{L2} we maximize the exterior region of the disk centered at $\lambda_i \in \Lambda_e$ (recall the shaded region in Fig.~\ref{fig1b}). If $0 \in \Lambda_i$, by the bounded-real lemma {\cite{Hespa}} \textbf{L1} should have at least one feasible solution provided $A$ in \eqref{eq:LTI_canonical} is Hurwitz.

We collect matrices $\Pi$ resulting from all feasible solutions to \textbf{L1} in the set $\mathbf{\Pi}_1(H_L)$, and those matrices resulting from all feasible solutions to \textbf{L2} in the set $\mathbf{\Pi}_2(H_L)$. Note that $(\mathbf{\Pi}_1(H_L)\cup \mathbf{\Pi}_2(H_L))\subset \mathbf{\Pi}(H_L)$, with $\mathbf{\Pi}(H_L)$ given in \eqref{eq:PI_set}. Through Theorem~\ref{th:general_SG_approx} this leads to the over-approximation
\begin{equation}
    \textup{SG}(H_L) \subseteq \bigcap_{\Pi \in (\mathbf{\Pi}_1(H_L)\cup \mathbf{\Pi}_2(H_L))}\mathcal{S}(\Pi),
\end{equation}
with the region $\mathcal{S}(\Pi)$ given in \eqref{eq:S_PI_canonical}.

In summary: an over-approximation of $\textup{SG}(H_L)$ with $H_L$ the LTI system in \eqref{eq:LTI_canonical} is obtained by solving both \textbf{L1} for $\lambda_c \in \Lambda_i$ and \textbf{L2} for $\lambda_c \in \Lambda_e$, and taking the intersection of regions resulting from all feasible solutions to these problems. We illustrate the procedure through the next example.

\begin{figure}[b]
		\centering
		\includegraphics[width=0.7\linewidth]{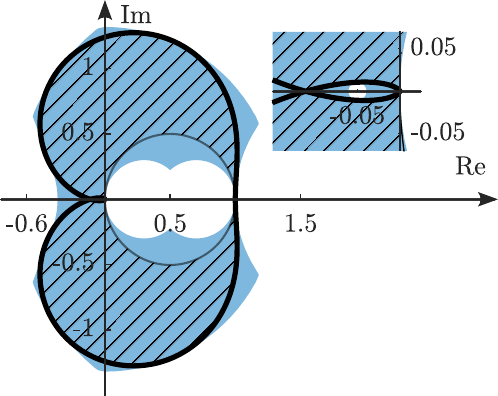}
		\caption{Nyquist diagram of $H_L(j\omega)$ in \eqref{eq:LTI_example} given by the thick black curve, $\textup{SG}(H_L)$ given by the hatched black region, and over-approximation using \eqref{eq:LiLeexample} given by the blue region.}
		\label{fig:crudeltiapprox}
	\end{figure}

\begin{example}\label{ex:LTI}
Consider a third-order LTI system $H_L$ described by the transfer function
\begin{align} \label{eq:LTI_example}
H_L(s) = \frac{1}{s^3+5s^2+2s+1}.
\end{align} 
Take 
\begin{subequations}\label{eq:LiLeexample}
\begin{align}
    \Lambda_i &= \{-1,0,1\},\\
        \Lambda_e  &= \{-2, -1, -0.05, 0.3,0.7,1,2  \}.
        \end{align}
        \end{subequations}
Solving \textbf{L1} and \textbf{L2} for the given sets in \eqref{eq:LiLeexample} leads to the region depicted in Fig.~\ref{fig:crudeltiapprox} in blue. From Theorem~\ref{th:LTI_bk_hull}, it follows that for SISO LTI systems, the exact scaled graph can be determined as the h-hull of the Nyquist diagram of $H_L(j\omega)$, $\omega \in \mathbb{R}$. Both the Nyquist diagram and $\textup{SG}(H_L)$ are show in Fig.~\ref{fig:crudeltiapprox} by the thick black line and the hatched black region, respectively. 
\end{example}

The quality of the over-approximation depends on the selection of the sets $\Lambda_i$ and $\Lambda_e$. In general, increasing the number of elements in these sets potentially yields a tighter over-approximation. To illustrate this point, we reconsider Example~\ref{ex:LTI} using different $\Lambda_i$ and $\Lambda_e$. 

\begin{figure}[b]
		\centering
		\includegraphics[width=0.7\linewidth]{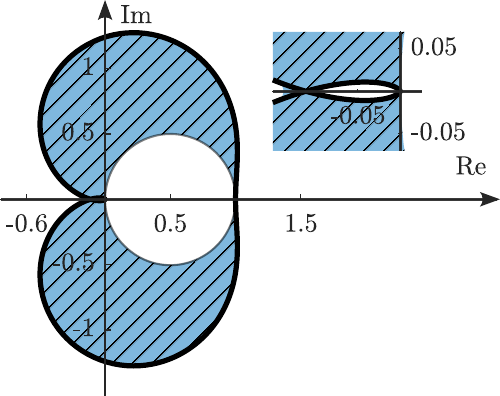}
		\caption{Nyquist diagram of $H_L$ in \eqref{eq:LTI_example} given by the thick black curve, $\textup{SG}(H_L)$ given by the hatched black region, and over-approximation using \eqref{eq:combined} given by the blue region.}
		\label{fig:ltisrgclose}
	\end{figure}

\noindent\textbf{Example 1 (revisited)}
\itshape Consider again the LTI system $H_L$ in \eqref{eq:LTI_example}, and select
\begin{subequations}\label{eq:combined}
\begin{align}
\Lambda_i &= \{-2+0.05k \mid k = 0,1, \dots, 80\},\\
{\Lambda}_e &=  \{-10+0.25k \mid k = 0,1, \dots, 80\}.
\end{align}
\end{subequations}
\noindent The over-approximation of $\textup{SG}(H_L)$ that results from solving \textbf{L1} and \textbf{L2} with \eqref{eq:combined} along with $\textup{SG}(H_L)$ are shown in Fig.~\ref{fig:ltisrgclose}. It can be clearly seen that, compared to the over-approximation in Fig.~\ref{fig:crudeltiapprox}, the over-approximation using the sets in \eqref{eq:combined} becomes significantly tighter, almost resembling $\textup{SG}(H_L)$ exactly.

\normalfont

 The improved tightness of the approximation in Example~\ref{ex:LTI} using larger sets $\Lambda_i$ and $\Lambda_e$ is not a mere consequence of the particular choice of the system $H_L$ in \eqref{eq:LTI_example}, but rather is a property of the connection between IQCs and SGs. In fact, it turns out that the approximation becomes exact for normal LTI systems, when setting $\Lambda_i=\mathbb{R}$ and $\Lambda_e =\mathbb{R}$. This will be formalized in the next section.

\subsection{Exactness}
To formally prove the previous statement on exactness, we require some additional concepts. 
\begin{defn}[half-planes in $\mathbb{R}^2$]\hfill
\begin{enumerate}
    \item  A Euclidean half-plane is an unbounded subset of $\mathbb{R}^2$ whose boundary is a straight line, and which partitions $\mathbb{R}^2$ into two disjoint regions.
    \item  A supporting half-plane $\mathcal{H}_s$ of a closed set $\mathcal{S} \subset \mathbb{R}^2$ is a Euclidean half-plane that fully contains $\mathcal{S}$, and whose boundary intersects $\mathcal{S}$ at least at one point.  
\end{enumerate}    
\end{defn}
We denote the set of all supporting half-planes of $\mathcal{S}\subset \mathbb{R}^2$ by ${\mathcal{H}(\mathcal{S})}$ and recall the following fundamental result.

\begin{lem}\label{lem:shy} A closed convex set $\mathcal{S}\subset \mathbb{R}^2$ having nonempty interior is the intersection of all its supporting half-planes, i.e., $\mathcal{S} = \cap_{\mathcal{H}_s \in \mathcal{H}(\mathcal{S})} \mathcal{H}_s$.
\end{lem}
Lemma~\ref{lem:shy} is a well-known result that follows from the converse supporting hyperplane theorem, see also {\cite[Section 2.5]{boyd2004convex}}. 

Equipped with the above definitions and results, we are ready to formulate the exactness result in this section.

\begin{figure}[htb]
	
	\centering
	
	\begin{subfigure}{0.49\linewidth}
	\centering\includegraphics[width=\linewidth]{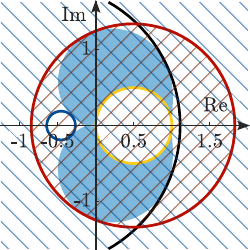}
	\caption{}
		\label{}
	\end{subfigure}
	\begin{subfigure}{0.49\linewidth}
	\centering\includegraphics[width=\linewidth]{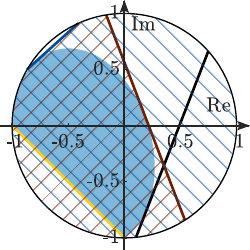}
	\caption{}
		\label{}
	\end{subfigure}
	\caption{(a) Illustration of $\textup{SG}(H)$ in blue, with several interior/exterior regions of a circle. (b) Image under the Beltrami-Klein mapping $f_{BK}$ \eqref{eq:bk_map}. Note that circles are mapped to straight line segments. The hatched regions indicate half-planes that contain (a) $\textup{SG}(H)$ and (b) $f_{BK}(\textup{SG}(H))$.}
		\label{fig:beltrami_map_example}
\end{figure}

\begin{thm}\label{th:final_LTI}
    Consider a normal LTI system $H_L$ of the form \eqref{eq:LTI_canonical} and assume that the matrix $A$ is Hurwitz and the pair $(A,B)$ is controllable. Let $\Lambda_i = \mathbb{R}$ and $\Lambda_e = \mathbb{R}$. Then
   \begin{equation}\label{eq:exact_LTI}
    {\overline{SG}}(H_L) = \bigcap_{\Pi \in (\mathbf{\Pi}_1(H_L)\cup \mathbf{\Pi}_2(H_L))}\mathcal{S}(\Pi),
\end{equation}
with $\mathbf{\Pi}_1(H_L), \mathbf{\Pi}_2(H_L)$ the sets containing all feasible solutions to problems \textbf{L1} and \textbf{L2}, respectively, and where $\mathcal{S}(\Pi)$ is given in \eqref{eq:S_PI_canonical}.
\end{thm}
\begin{pf}
    Let us denote $\overline{\textup{SG}}(H_L)_+ = \overline{\textup{SG}}(H_L)\cap \mathbb{C}_+$. First, note that by virtue of Theorem \ref{th:LTI_bk_hull} we have that $f_{BK}(\overline{\textup{SG}}(H_L)_+)$ is a convex set that is contained in the unit disk $\mathbb{D}$. It then follows from Lemma~\ref{lem:shy} that 
\begin{equation}
   {f_{BK}}(\overline{\textup{SG}}(H_L)_+) = \bigcap_{\mathcal{H}_s \in \mathcal{H}({f_{BK}}(\overline{\textup{SG}}(H_L)_+))} (\mathcal{H}_s \cap \mathbb{D}).\end{equation}
    Since $f_{BK}^{-1}(f_{BK}(z)) = \left\{z,\bar{z}\right\}$ for $z \in \mathbb{C}$ and the fact that scaled graphs are symmetric around the real axis, we find
  \begin{equation}\label{eq:34}
   \overline{\textup{SG}}(H_L) = {f^{-1}_{BK}}\left(\bigcap_{\mathcal{H}_s \in \mathcal{H}({f_{BK}}(\overline{\textup{SG}}(H_L)))}(\mathcal{H}_s \cap \mathbb{D})\right).\end{equation}  
Since for $a,b \in \mathbb{D}$, $f^{-1}_{BK}(a) \cap f^{-1}_{BK}(b) \neq \emptyset$ implies $a = b$ (i.e., $f^{-1}_{BK}$ is set-valued injective on $\mathbb{D}$), it follows that for $A,B \subset  \mathbb{D}$, $f^{-1}_{BK}(A\cap B) = f^{-1}_{BK}(A)\cap f^{-1}_{BK}(B)$. Hence, we conclude
\begin{equation}\label{eq:36}
      \overline{\textup{SG}}(H_L) = \bigcap_{\mathcal{H}_s \in \mathcal{H}({f_{BK}}(\overline{\textup{SG}}(H_L)))}{f_{BK}^{-1}}(\mathcal{H}_s \cap \mathbb{D}).
\end{equation}
When $\mathcal{H}_s\cap \mathbb{D}= \mathbb{D}$, then $f_{BK}^{-1}(\mathcal{H}_s\cap \mathbb{D}) = \mathbb{C}$. On the other hand, when $\mathcal{H}_s \cap \mathbb{D} \neq \mathbb{D}$, the boundary of this region is composed by part of the unit circle and a straight line segment (as a result of the intersection of a half-plane and a unit disk). Since a straight line segment in $\mathbb{D}$ maps via $f^{-1}_{BK}$ to a circle in $\mathbb{C}$ centered on the real axis, it follows that, in this case, $f^{-1}_{BK}(\mathcal{H}_s\cap \mathbb{D})$ maps to either the interior or exterior of such circle. In other words, $f^{-1}_{BK}(\mathcal{H}_s\cap \mathbb{D})= \mathcal{S}(\Pi(\sigma, \lambda_c, r))$, with $\Pi$ as in \eqref{eq:PI_parametrization}, $\sigma \in \left\{-1,+1\right\}$, $\lambda_c \in \mathbb{R}$, and $r\geq 0$. Collect all matrices $\Pi$ corresponding to these regions in $\mathbf{\hat{\Pi}}(H_L)$. Then, we find
\begin{equation}
    \overline{\textup{SG}}(H_L) = \bigcap_{\Pi \in \mathbf{\hat{\Pi}}(H_L)}\mathcal{S}(\Pi).
\end{equation}
By Lemma~\ref{lem:connection_to_SG}, each $\mathcal{S}(\Pi)$ corresponds to an IQC of the form \eqref{eq:Supply_inequality}, which, by Theorem~\ref{th:LTI_approx} implies feasibility of the LMI in \eqref{eq:LMI_F}. Since $\Lambda_i = \mathbb{R}$ and $\Lambda_e = \mathbb{R}$, and since we minimize/maximize $r$ in \textbf{L1} and \textbf{L2}, it follows that $\mathbf{\hat{\Pi}}(H_L) \subset (\mathbf{\Pi}_1(H_L)\cup \mathbf{\Pi}_2(H_L))$. Note that all matrices $\Pi$ that are not contained in $\mathbf{\hat{\Pi}}(H_L)$ but are contained in $\mathbf{\Pi}_1(H_L)\cup \mathbf{\Pi}_2(H_L)$ necessarily provide an excessive over-approximation of $\textup{SG}(H_L)$. Since we take the intersection of regions, we find that
\begin{equation}
    \overline{\textup{SG}}(H_L)=\!\!\bigcap_{\Pi \in \mathbf{\hat{\Pi}}(H_L)}\!\!\mathcal{S}(\Pi) = \!\!\!\!\bigcap_{\Pi \in (\mathbf{\Pi}_1(H_L)\cup \mathbf{\Pi}_2(H_L))} \!\!\! \mathcal{S}(\Pi)
\end{equation}
and thus \eqref{eq:exact_LTI}. \hfill $\blacksquare$
\end{pf}
Explicit mappings from interior/exterior regions of a circle in $\mathbb{C}$ to the unit disk $\mathbb{D}$ and vice versa are shown in Fig.~\ref{fig:beltrami_map_example} to highlight the main mechanisms used in the proof of Theorem~\ref{th:final_LTI}. 

\begin{rem}
    In Theorem~\ref{th:final_LTI} we work with the closure of the scaled graph. Although in some cases $\textup{SG}(H)=\overline{\textup{SG}}(H)$ (see, e.g., the example in the remark below), the scaled graph might not be a closed set in general. Despite taking the closure, Theorem~\ref{th:final_LTI} still shows that we can construct an arbitrarily tight over-approximation of $\textup{SG}(H)$ for normal LTI systems. 
\end{rem}

\begin{rem}\label{rem:1D}
Note that, in general, problem \textbf{L1} and \textbf{L2} may require solving an infinite number of LMIs for obtaining the exact scaled graph of a system. However, for certain systems an exact representation can be obtained by solving a finite number of LMIs. An example of such systems is given by first-order LTI systems having a transfer function of the form
    \begin{align} 
    H_L(s) = \frac{1}{s+a},\quad \text{ with } a\in\mathbb{R}_{>0}.
\end{align} 
In this case, the scaled graph is given by
\begin{equation*}
\begin{split}&\textup{SG}(H_L)=\\
&\left\{z\in \mathbb{C} \mid \left(\textup{Re}\left\{z\right\}-\frac{1}{2a}\right)^2+\textup{Im}\left\{z\right\}^2-\frac{1}{4a^2} = 0 \right\},\end{split} 
\end{equation*}
i.e., a shifted circle with center $\lambda_c = \tfrac{1}{2a}$ and radius $r = \tfrac{1}{2a}$. As such, the exact scaled graph can be obtained by solving \textbf{L1} and \textbf{L2} with $\Lambda_i = \Lambda_e =\left\{\tfrac{1}{2a}\right\}$ leading to the solutions $P = \tfrac{1}{2a}$, $r = \tfrac{1}{2a}$ for \textbf{L1}, and $P=-\tfrac{1}{2a}$, $r = \tfrac{1}{2a}$ for \textbf{L2}.
\end{rem}

The result in Theorem \ref{th:final_LTI} shows that the computation procedure proposed in Section~\ref{sec:proc} is arbitrary tight for normal, stable LTI systems. This result gives a solid foundation for extending the approach towards computation of soft scaled graphs for classes of nonlinear systems as well. We exploit this in the next section.

\section{Results for Classes of Nonlinear Systems}
\label{sec:NL}
In this section, we extend the methods introduced in the previous section to compute (over-approximations of) the scaled graph for certain relevant classes of nonlinear systems. We start with reset/impulsive systems, and subsequently consider the class of piecewise linear systems. 

\subsection{Reset systems}
    Consider a reset system of the form \cite{aangenent_performance_2010}, \cite{zaccarian_first_2005}, \cite{ZHANG2024106063}
\begin{align}
	H_R: \left\{\begin{aligned}
		\dot{x}(t)  &= Ax(t)+Bu(t), & \textnormal{if}\; \xi(t)\in\mathcal{F},\\
		x(t^+) &=Rx(t), & \textnormal{if}\;\xi(t)\in \mathcal{J},\\
		y(t) &=Cx(t)+Du(t), &
	\end{aligned} \right. \label{eq:canonical_reset_system}
\end{align}
with states $x(t)\in\mathbb{R}^m$, initial condition $x(0)=0$, input $u(t)\in\mathbb{R}^n$, output $y(t)\in \mathbb{R}^n$, and $\xi(t) = [x^\top (t), u{^\top}(t)]^\top \in \mathbb{R}^{m+n}$ all at time $t \in \mathbb{R}_{\geq 0}$, and where $x(t^+) = \lim_{s \downarrow t}x(s)$. The matrix $R \in \mathbb{R}^{m\times m}$ defines the reset map, and the jump and flow sets are given by
\begin{align}
	\mathcal{F} &= \{\mu \in \mathbb{R}^{m+n} \;|\; \mu^TM\mu\geq 0\} \label{eq:canonical_flow_set_reset},\\
	\mathcal{J} &= \{\mu \in \mathbb{R}^{m+n} \;|\; \mu^TM\mu\leq 0\} \label{eq:canonical_jump_set_reset},
\end{align}
with $M\in\mathbb{S}^{m+n}$. Furthermore, we associate with $H_R$ in \eqref{eq:canonical_reset_system} a base linear system given by $(A,B,C,D)$, i.e., an LTI system as in \eqref{eq:LTI_canonical} without resets. We assume that solutions to \eqref{eq:canonical_reset_system} are well defined for all times $t\in\mathbb{R}_{\geq0}$, and resets occur at isolated points in time. Note that the above assumptions may only hold true for specific inputs that satisfy certain regularity conditions, see, e.g., the discussion in \cite[Section 3.2]{aangenent_performance_2010}. Although for LTI systems this is not an issue, to obtain a meaningful scaled graph of the reset system in \eqref{eq:canonical_reset_system}, we may be restricted to consider a subset $\mathcal{U}$ of the input space $\mathcal{L}_2$. Furthermore, we assume that $u \in \mathcal{U}$ leads to $x,y\in\mathcal{L}_2$, which can be verified through several techniques such as in \cite{aangenent_performance_2010}, \cite{zaccarian_first_2005}.

\begin{thm}
	\label{th:improved_reset_KYP}
	Consider a reset system $H_R:\mathcal{L}_2\rightrightarrows\mathcal{L}_2$ of the form \eqref{eq:canonical_reset_system}. Suppose there exist $\tau_1,\tau_2\geq 0$ and a matrix $P\in \mathbb{S}^n$ that satisfy the LMIs
    \begin{subequations}
        \begin{align}
         \begin{bmatrix}
			A &B\\I&0
		\end{bmatrix}^\top\begin{bmatrix}
			0 & P \\ P &0
		\end{bmatrix}\begin{bmatrix}
			A &B\\I&0
		\end{bmatrix}-\Theta(\Pi)+\tau_1  M&\preceq 0, \label{eq:lmi_reset1}\\
        \begin{bmatrix}
            R^\top PR-P & 0 \\
            0 & 0
        \end{bmatrix}-\tau_2 M & \preceq 0,\label{eq:lmi_reset2}
        \end{align}
    \end{subequations}
    with $\Theta(\Pi)$ given as in \eqref{eq:thetaPi}. Then, $\textup{SG}_{\mathcal{U}}(H_R) \subseteq \mathcal{S}(\Pi)$, with $\mathcal{S}(\Pi)$ as in \eqref{eq:S_PI_canonical}. 
\end{thm}
\begin{pf}
Suppose the conditions hypothesized in the theorem hold true. Consider a quadratic function of the form $W(x(t)) = x(t)^\top Px(t)$. Pre- and post-multiplying \eqref{eq:lmi_reset1} with $\xi(t)^\top = [x(t)^\top, u(t)^\top]$ and $\xi(t)$, respectively, implies for almost all $t$ during flow ($\xi(t)\in \mathcal{F}$)
\begin{equation}\label{eq:dWdt}
    \frac{d}{dt}W(x(t))\leq \xi(t)^\top \Theta(\Pi)\xi(t),
\end{equation}
where we have used the fact that $\tau_1\xi(t)^\top M \xi(t) \geq 0$ when $\xi(t) \in \mathcal{F}$. In a similar manner, the inequality in \eqref{eq:lmi_reset2} implies
\begin{equation}\label{eq:Wreset}
    W(x(t^+))\leq W(x(t)),
\end{equation}
when $\xi(t) \in \mathcal{J}$. Integrating \eqref{eq:dWdt} over a time interval that the reset system is flowing, and combining this with non-increasing of $W$ along resets as in \eqref{eq:Wreset} and $x(0) = 0$, yields
\begin{equation}\label{eq:Wint}
    W(x(T)) \leq \int_{0}^T\xi(t)^\top \Theta(\Pi)\xi(t)dt.
\end{equation}
Under the assumption that resets occur at isolated times, the signal $x$ is differentiable almost everywhere. Furthermore, since by assumption $u \in \mathcal{U}\subset \mathcal{L}_2$ implies $x \in \mathcal{L}_2$, it follows that $\dot{x} \in \mathcal{L}_2$. Invoking \cite[p. 237]{Desoer} shows that $\lim_{t\to \infty}x(t)=0$. Taking $T \to \infty$ in \eqref{eq:Wint} leads to $\lim_{T\to \infty}W(x(T))=0$, and, therefore,
\begin{equation}
    0 \leq \int_{0}^\infty\begin{bmatrix} y(t) \\ u(t) \end{bmatrix}^\top \left(\Pi \otimes I_n \right)\begin{bmatrix} y(t) \\ u(t) \end{bmatrix}dt.
\end{equation}
Applying Lemma~\ref{lem:connection_to_SG} leads to the result. \hfill $\blacksquare$
\end{pf}
Denote by $F(\sigma,\lambda_c,r,P,\tau_1)$ the matrix in the left-hand side of the inequality in \eqref{eq:lmi_reset1}, and by $G(P,\tau_2)$ the matrix in the left-hand side of the inequality in \eqref{eq:lmi_reset2}. Theorem~\ref{th:improved_reset_KYP} allows to formulate a soft scaled graph computation method akin to the one for LTI systems in Section~\ref{sec:proc} as the following optimization problems: 

\textbf{R1:} For each $\lambda_c \in \Lambda_i$, 
\begin{equation}\label{eq:optProbreset}
\begin{aligned}
&\underset{r,P, \tau_1, \tau_2}{\text{minimize}}        & & r^2\\
&\text{subject to} & &F(-1,\lambda_c,r,P,\tau_1)\preceq 0, \\
& & & G(P,\tau_2) \preceq 0,\\
& & & P \in \mathbb{S}^n, \tau_1,\tau_2,r> 0.
\end{aligned}
\end{equation}
\textbf{R2:} For each $\lambda_c \in \Lambda_e$, 
\begin{equation}\label{eq:optProb}
\begin{aligned}
&\underset{r,P,\tau_1,\tau_2}{\text{maximize}}        & & r^2\\
&\text{subject to} & &F(+1,\lambda_c,r,P,\tau_1)\preceq 0, \\
& & & G(P,\tau_2) \preceq 0,\\
& & & P \in \mathbb{S}^n, \tau_1,\tau_2,r> 0.
\end{aligned}
\end{equation}

Again, we collect matrices $\Pi$ resulting from all feasible solutions to \textbf{R1} in the set $\mathbf{\Pi}_1(H_R)$, and those matrices resulting from all feasible solutions to \textbf{R2} in the set $\mathbf{\Pi}_2(H_R)$, such that to over-approximation is given by
\begin{equation}
    \textup{SG}(H_R) \subseteq \bigcap_{\Pi \in (\mathbf{\Pi}_1(H_R)\cup \mathbf{\Pi}_2(H_R))}\mathcal{S}(\Pi),
\end{equation}
with the region $\mathcal{S}(\Pi)$ given in \eqref{eq:S_PI_canonical}.

To demonstrate effectiveness of the computation method also in the reset case, consider the next example.

\begin{example}\label{ex:reset}
We adopt the reset system from \cite[Example 1]{van_den_eijnden_scaled_2024} described by \eqref{eq:canonical_reset_system} with the data
\begin{align}
		\left[\begin{array}{c|c }
			A&B   \\
			\hline 
			C&D		
		\end{array}\right] = 
		\left[\begin{array}{c c|c }
			 -1 & 0 & 1   \\
			 1 & -1 & 0\\
			\hline 
			0 & 1 & 0
		\end{array}\right],
			\end{align}
			and $R=\textup{diag}(0,0)$, $M=\textup{diag}(0.9^2,-1,0).$
    
Solving \textbf{R1} and \textbf{R2} with 
\begin{align}
\Lambda_i &= \{-1+0.05k \mid k = 0,1, \dots, 80\},\\
{\Lambda}_e &=  \{-1+0.25k \mid k = 0,1, \dots, 80\},
\end{align}
results in the over-approximation of $\textup{SG}(H_R)$ depicted in Fig.~\ref{fig:resetcomparrison} by the blue region. The over-approximation obtained through the method in \cite{van_den_eijnden_scaled_2024} is indicated by the hatched black region. Clearly, the new approach leads to a significantly reduced over-approximation, thereby highlighting the effectiveness.  
	\begin{figure}[htb]
		\centering
		\includegraphics[width=0.5\linewidth]{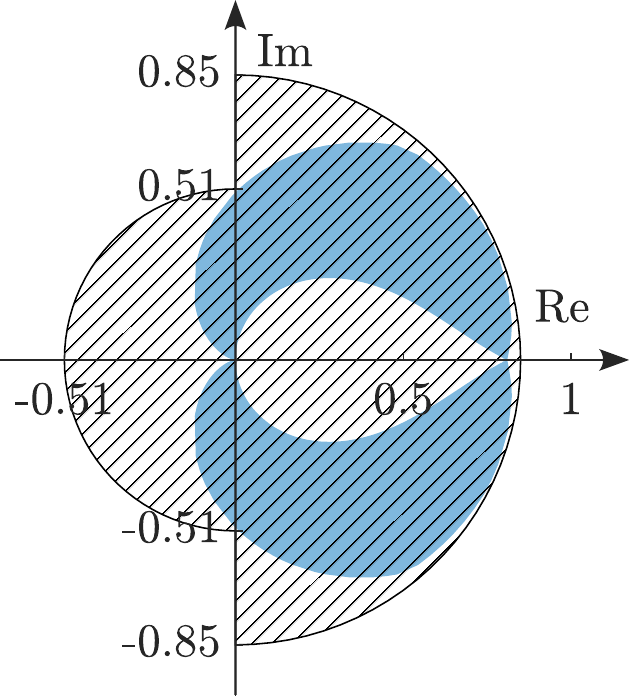}
		\caption{Over-approximation of $\textup{SG}(H_R)$ with $H_R$ the reset system in \eqref{eq:canonical_reset_system} using i) the method in \cite{van_den_eijnden_scaled_2024} (hatched black region) and ii) the new procedure (blue region).}
		\label{fig:resetcomparrison}
	\end{figure}

To gain further insights in tightness of the over-approximation, an under-approximation of $\textup{SG}(H_R)$ that results from numerically ``sampling'' the scaled graph with multi-sine inputs of the form
\begin{align} \label{eq:sample_inputs}
    u(t) = \left(\sum_{m=1}^{M}k_m\sin(\omega_mt+\varphi_m)\right)e^{\mu t},
    \end{align}
with $M\leq 20$, $k_m,\omega_m,\varphi_m$, and $\mu<0$ randomized parameters, is shown in Fig.~\ref{fig:reset_a} by the collection of blue dots. A clear discrepancy is visible between the under- and over-approximation. Such discrepancy is expected from i) the fact that our approach only guarantees an over-approximation for reset systems (the real SG may be smaller), and ii) the limited nature of the input signals in \eqref{eq:sample_inputs} (the real SG will be larger than the region indicated by the blue dots). 

To indicate the effects of adding resets to LTI systems, Fig.~\ref{fig:reset_b} shows the scaled graph of the base linear system described by the matrices $(A,B,C,D)$. Interestingly, the reset action deforms the scaled graph of the base linear system, altering its gain-phase properties. Note that the maximum gain of $1$ of the base linear system is reduced as a consequence of the reset rule. 

\begin{figure}[htb]
	
	\centering
	
	\begin{subfigure}{0.49\linewidth}
	\includegraphics[width=\linewidth]{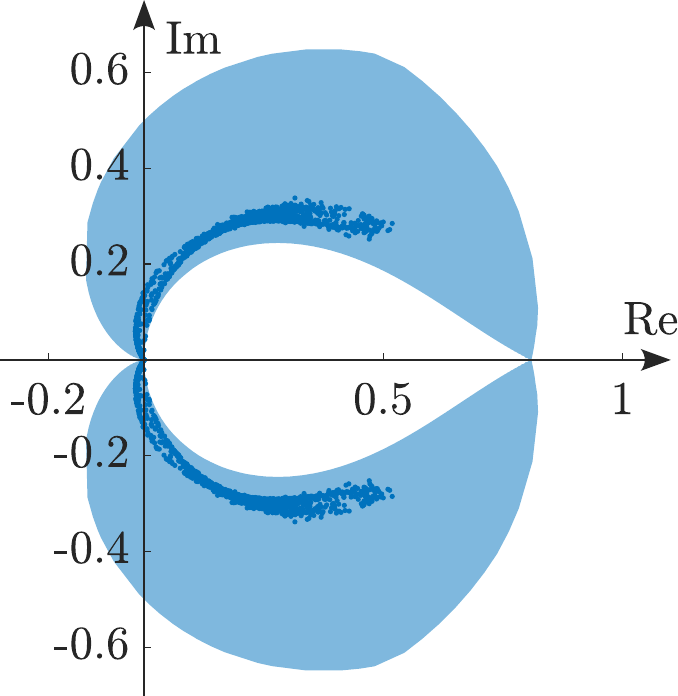}
	\caption{}
    \label{fig:reset_a}
	\end{subfigure}
	\begin{subfigure}{0.49\linewidth}
	\includegraphics[width=\linewidth]{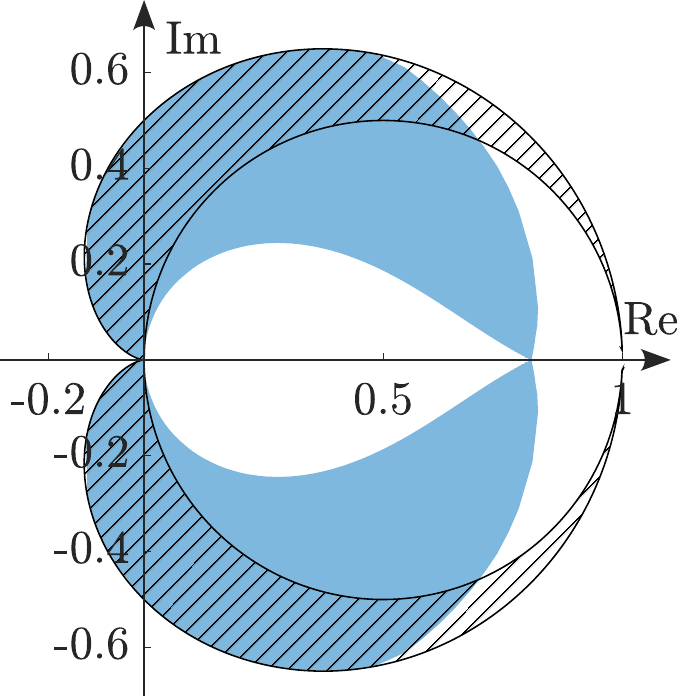}
	\caption{}
    \label{fig:reset_b}
	\end{subfigure}
	\caption{(a) Over-approximation of $\textup{SG}(H_R)$ with $H_R$ the reset system in \eqref{eq:canonical_reset_system} (blue region), and numerically sampled under-approximation (blue dots). (b) Scaled graph of the base linear system (hatched black region), along with the over-approximation of $\textup{SG}(H_R)$ (in blue).}
		\label{fig:resetsampling}
\end{figure}
   \end{example}

\subsection{Piecewise linear systems}
Next, we consider piecewise linear (PWL) systems of the form
\begin{equation}\label{eq:PWL}
    H_P: \begin{cases}
        \dot{x}(t) & = A_i x(t) +B_i u(t), \\
        y(t)& = C_i x(t)+ D_i u(t),
    \end{cases}\textup{ if }\xi(t) \in \mathcal{C}_i, \:i \in \mathcal{N},
\end{equation}
with states $x(t) \in \mathbb{R}^m$, initial condition $x(0) = 0$, input $u(t)\in\mathbb{R}^n$, output $y(t)\in \mathbb{R}^n$, and $\xi(t) = [x^\top (t), u(t)]^\top \in \mathbb{R}^{m+n}$ all at time $t \in \mathbb{R}_{\geq 0}$. The flow sets are given by
\begin{equation}
    \mathcal{C}_i = \left\{\mu \in \mathbb{R}^{m+n} \mid E_i \mu \geq 0\right\},
\end{equation}
where $E_i \in \mathbb{R}^{p\times (m+n)}$, $i\in \mathcal{N}$, and $\mathcal{N}=\left\{1,\ldots,N\right\}$. We assume that $\cup_{i \in \mathcal{N}}\;\mathcal{C}_i = \mathbb{R}^{m+n}$, and that solutions to \eqref{eq:PWL} exist in the sense of Carath\'eodory, i.e., absolutely continuous functions that satisfy \eqref{eq:PWL} for almost all times. Moreover, it is assumed that $u \in \mathcal{U}\subset\mathcal{L}_2$ leads to $x,y \in \mathcal{L}_2$, such that $\dot{x} \in \mathcal{L}_2$. 
\begin{thm}
	\label{th:improved_PWL_KYP}
	Consider a PWL system $H_P:\mathcal{L}_2 \rightrightarrows \mathcal{L}_2$ of the form \eqref{eq:PWL}. Suppose there exist matrices $U_i\in \mathbb{S}_{\geq 0}^p$, $i\in \mathcal{N}$, and a matrix $P\in \mathbb{S}^n$ that satisfy the LMIs
        \begin{align}
         \begin{bmatrix}
			A_i &B_i\\I&0
		\end{bmatrix}^\top\begin{bmatrix}
			0 & P \\ P &0
		\end{bmatrix}\begin{bmatrix}
			A_i &B_i\\I&0
		\end{bmatrix}-\Theta_i(\Pi)+E_i^\top U_iE_i&\preceq 0, \label{eq:lmi_PWL}
        \end{align}
 for all $i \in \mathcal{N}$, with
 \begin{equation}
    \Theta_i(\Pi) = \begin{bmatrix}
		C_i & D_i \\ 0 & I
	\end{bmatrix}^\top(\Pi\otimes I_n)\begin{bmatrix}
		C_i & D_i \\ 0 & I
	\end{bmatrix}, 
 \end{equation}
 and $\Pi \in \mathbb{S}^2$. Then, $\textup{SG}_\mathcal{U}(H_P) \subseteq \mathcal{S}(\Pi)$, with $\mathcal{S}(\Pi)$ in \eqref{eq:S_PI_canonical}. 
\end{thm}
\begin{pf}
Suppose the conditions hypothesized in the theorem hold true. Consider a quadratic function of the form $W(x(t)) = x(t)^\top Px(t)$. Pre- and post-multiplying \eqref{eq:lmi_PWL} with $\xi(t)^\top = [x(t)^\top, u(t)^\top]$ implies for almost all time $t$ when $\xi(t) \in \mathcal{C}_i$
\begin{equation}\label{eq:dWdt2}
    \frac{d}{dt}W(x(t))\leq \xi(t)^\top \Theta_i(\Pi)\xi(t), \: i\in \mathcal{N},
\end{equation}
where we have used the fact that since $U_i$ has nonnegative elements, $\xi(t)^\top E_i^\top U_iE_i \xi(t) \geq 0$ when $\xi(t) \in \mathcal{C}_i$, $i\in \mathcal{N}$. Integrating \eqref{eq:dWdt2} over a time interval $[0,T]$, and using $x(0) = 0$ such that $W(x(0)) = 0$ yields
\begin{equation}\label{eq:Wint2}
    W(x(T)) \leq \int_{0}^T\begin{bmatrix} y(t) \\ u(t) \end{bmatrix}^\top \left(\Pi \otimes I_n \right)\begin{bmatrix} y(t) \\ u(t) \end{bmatrix}dt.
\end{equation}
Since by assumption $x \in \mathcal{L}_2$, and $\dot{x} \in \mathcal{L}_2$ almost everywhere, we have by virtue of \cite[p. 237]{Desoer} that $\lim_{t\to \infty}x(t)=0$. Taking $T \to \infty$ in \eqref{eq:Wint2} leads to $\lim_{T\to \infty}W(x(T))=0$, and, therefore,
\begin{equation}
    0 \leq \int_{0}^\infty\begin{bmatrix} y(t) \\ u(t) \end{bmatrix}^\top \left(\Pi \otimes I_n \right)\begin{bmatrix} y(t) \\ u(t) \end{bmatrix}dt.
\end{equation}
The result follows from applying Lemma~\ref{lem:connection_to_SG}.\hfill $\blacksquare$\end{pf}
As before, based on Theorem~\ref{th:improved_PWL_KYP} a procedure for over-approximating the scaled graph of the PWL system in \eqref{eq:PWL} can be formulated by reformulating Problem 1 and Problem 2 in Section~\ref{sec:proc} accordingly. Denote by $F_i(\sigma,\lambda_c,r,P,U_i)$ the matrix in the left-hand side of \eqref{eq:lmi_PWL}. The procedure is formulated as follows.

\textbf{P1:} For each $\lambda_c \in \Lambda_i$, 
\begin{equation}\label{eq:optProbreset1}
\begin{aligned}
&\underset{r,P, U_i}{\text{minimize}}        & & r^2\\
&\text{subject to} & &F_i(-1,\lambda_c,r,P,U_i)\preceq 0, \:i \in \mathcal{N}\\
& & & P \in \mathbb{S}^n, U_i\in \mathbb{S}_{\geq 0}^p,r\geq 0.
\end{aligned}
\end{equation}
\textbf{P2:} For each $\lambda_c \in \Lambda_e$, 
\begin{equation}\label{eq:optProb2}
\begin{aligned}
&\underset{r,P,U_i}{\text{maximize}}        & & r^2\\
&\text{subject to} & &F_i(+1,\lambda_c,r,P,U_i)\preceq 0, \:i \in \mathcal{N},\\
& & & P \in \mathbb{S}^n, U_i\in \mathbb{S}_{\geq 0}^p,r\geq 0.
\end{aligned}
\end{equation}
Collecting all feasible solutions to the above problems allows to construct an over-approximation of the PWL system in \eqref{eq:PWL}, as we demonstrate in the next example. 
\begin{example} \label{ex:PWL}
Consider the PWL system $H_P$, with $\mathcal{N}=\{1,2,3,4\}$ and data
\begin{align}
		\left[\begin{array}{c|c }
			A_1&B_1   \\
			\hline 
			C_1&D_1		
		\end{array}\right]& =\left[\begin{array}{c|c }
			A_3&B_3   \\
			\hline 
			C_3&D_3		
		\end{array}\right] = 
		\left[\begin{array}{c c|c }
			 -0.1 & 0 & 0   \\
			 -1 & -2 & 1\\
			\hline 
			1 &0 & 0
		\end{array}\right], \\
         \left[\begin{array}{c|c }
			A_2&B_2   \\
			\hline 
			C_2&D_2		
		\end{array}\right] &=\left[\begin{array}{c|c }
			A_4&B_4   \\
			\hline 
			C_4&D_4		
		\end{array}\right] = 
		\left[\begin{array}{c c|c }
			 -0.1 & 1 & 0   \\
			 -1 & -1 & 1\\
			\hline 
			1 &0 & 0
		\end{array}\right],
\label{eq:pwl_parameters}   
	\end{align}
and
\begin{align}
E_1 &= -E_3 = \begin{bmatrix}1 & 0 & 0\\ 0&1 & 0\end{bmatrix} , \; E_2 = -E_4 = \begin{bmatrix}-1 & 0 & 0\\ 0&1 & 0 \end{bmatrix}.
\end{align} 
Solving \textbf{P1} and \textbf{P2} with  
\begin{align} 
\Lambda_i &= \{-20+0.1k \mid k = 0,1, \dots, 210\},\\
{\Lambda}_e &=  \{-50+0.35k \mid k = 0,1, \dots, 300\},
\end{align} 
results in the over-approximation of $\textup{SG}(H_P)$ shown in Fig~\ref{fig:a} by the region in blue. To obtain insights in tightness, again an under-approximation of the scaled graph that results from sampling the input-output behaviour using multi-sines of the form in \eqref{eq:sample_inputs} is constructed, as shown in Fig.~\ref{fig:a} by the blue dots. The discrepancy seen in this figure comes from potential conservatism, and the limited number of sampled input signals. 

\begin{figure}[hbt]
	
	\centering
	
	\begin{subfigure}{0.45\linewidth}
    \centering
	\includegraphics[width=\linewidth]{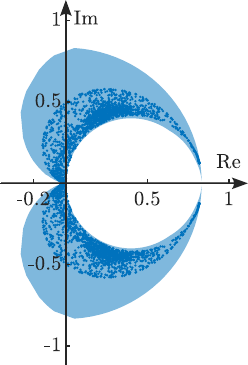}
	\caption{}
    \label{fig:a}
	\end{subfigure}
	\begin{subfigure}{0.45\linewidth}
    \centering
	\includegraphics[width=\linewidth]{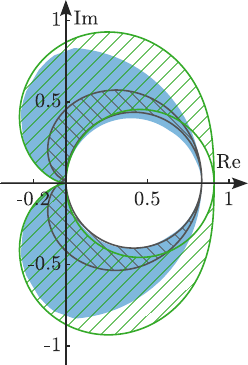}
	\caption{}
    \label{fig:b}
	\end{subfigure}
	\caption{(a) Over-approximation of $\textup{SG}(H_P)$ with $H_P$ the PWL system in \eqref{eq:PWL} (blue region), and numerically sampled under-approximation (blue dots). (b) Scaled graphs of mode $i=1$ (hatched black region), and mode $i=2$ (hatched green region), and over-approximation of $\textup{SG}(H_P)$ (in blue).}
\label{fig:PWL_example_sampled}
\end{figure}

To gain additional insights in the effects of switching, the scaled graphs of the individual LTI subsystems described by the matrices $(A_1, B_1,C_1,D_1)$ and $(A_2,B_2,C_2,D_2)$ are depicted in Fig.~\ref{fig:b} by the hatched black and hatched green regions, respectively. It is interesting to observe that, under the given switching law, the scaled graph of the switched system appears as a clear combination of the graphs of the two subsystems. 
\end{example}

\begin{figure*}[t]
	
	\centering
	
	\begin{subfigure}{0.31\linewidth}
    \centering
	\includegraphics[width=0.8\linewidth]{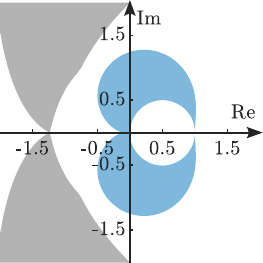}
	\caption{$\textup{SG}^{\dagger}(-H_R)$ (grey), $\textup{SG}(H_L)$ (blue)}
    \label{fig:FB_PWL_reset}
	\end{subfigure}
	\begin{subfigure}{0.31\linewidth}
    \centering
	\includegraphics[width=.8\linewidth]{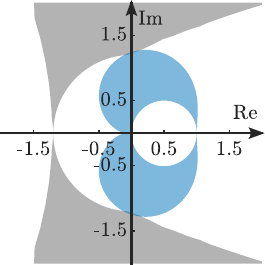}
	\caption{$\textup{SG}^{\dagger}(-H_P)$ (grey), $\textup{SG}(H_L)$ (blue)}
    \label{fig:FB_PWL_LTI}
	\end{subfigure}
    \begin{subfigure}{0.31\linewidth}
    \centering
	\includegraphics[width=.8\linewidth]{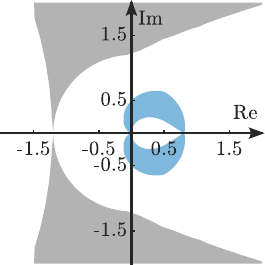}
	\caption{$-\textup{SG}^{\dagger}(H_P)$ (grey), $\textup{SG}(H_R)$ (blue)}
    \label{fig:FB_reset_LTI}
	\end{subfigure}
    
	\caption{Verification of Theorem~\ref{th:SGfb} for various feedback interconnections of the systems considered in Examples~\ref{ex:LTI}, \ref{ex:reset}, and \ref{ex:PWL}. (a) LTI and reset systems $(H_1,H_2)=(H_R,H_L)$, (b) LTI and PWL systems $(H_1,H_2)=(H_P,H_L)$, (c) reset and PWL systems $(H_1,H_2)=(H_P,H_R)$. In all figures, the grey regions depict $\textup{SG}^\dagger(-H_1)$ and the blue regions depict $\textup{SG}(H_2)$.}
\label{fig:interconnection_results}
\end{figure*}
We remark that the LMI-based computational procedures outlined above are not limited to te systems considered, but can be adapted to any class of systems having an LTI-like structure, including discrete-time LTI systems, linear-parameter varying systems \cite{toth2010modeling}, and Lur'e systems \cite{degroot2025exploitingstructuremimoscaled}. 

\subsection{Interconnection results}\label{sec:inter}
To illustrate the use of the previously computed scaled graphs for feedback analysis, we study stability on the basis of Theorem~\ref{th:SGfb} for three feedback interconnections\footnote{Here, we separate the question for stability from the question of (global) existence of solutions; the latter is beyond the scope of the current paper and well-posedness is assumed.} of the form as depicted in Fig.~\ref{fig:FB}: 
\begin{enumerate}[label=(\alph*)]
    \item $(H_1,H_2)=(H_R,H_L)$ with $H_L$ and $H_R$ the LTI and reset systems from Examples~\ref{ex:LTI} and~\ref{ex:reset};
    \item $(H_1,H_2)=(H_P,H_L)$ with $H_L$ and $H_P$ the LTI and PWL systems from Examples~\ref{ex:LTI} and~\ref{ex:PWL};
    \item $(H_1,H_2)=(H_P,H_R)$ with $H_R$ and $H_P$ the reset and PWL systems from Examples~\ref{ex:reset} and~\ref{ex:PWL}.
\end{enumerate}
The corresponding results from verifying Theorem~\ref{th:SGfb} with the previously obtained SG over-approximations for each example are shown in Fig.~\ref{fig:interconnection_results}.
It can be seen that stability is verified for case (a) and case (c), and the respective $\mathcal{L}_2$-gains are given by $1/r = 1/0.3=2.9$ and $1/r=1/0.66=1.5$. Note that the interconnection in case (c) seems significantly more robust with respect to gain variations compared to case (a). Stability for case (b) cannot be verified through Theorem~\ref{th:SGfb} as $\textup{SG}^\dagger(-H_P)$ and $\textup{SG}(H_L)$ overlap. However, the graphical nature of scaled graphs gives insight into possible directions for (re)design of the system/controller to achieve (robust) stability. For example, $H_L$ might be scaled by a proportional gain to shrink $\textup{SG}(H_L)$ or combined with an LTI lead filter to rotate $\textup{SG}(H_L)$ towards the right-half plane, see also \cite{chaffey_loop_2022}. We stress that the above performance, robustness, and design insights follow easily from the plots in Fig.~\ref{fig:interconnection_results}, but are generally difficult to obtain from a pure LMI-based stability analysis. Besides, the above examples demonstrate how scaled graphs allow for easy graphical stability verification for interconnections of different classes of systems. 

\section{Computations for hard scaled graphs}\label{sec:hardSG}
So far, we have focused on computing soft scaled graphs; recall Section~\ref{sec:sg}. In this section, we briefly comment on adapting the foregoing results for the computation of hard scaled graphs.

\subsection{Hard scaled graphs and hard IQCs}
The hard scaled graph of a system $H:\mathcal{L}_{2e}\rightrightarrows\mathcal{L}_{2e}$ over inputs $u \in {\mathcal{U}_e}\subset\mathcal{L}_{2e}$ is given by \cite{chen2025softhardscaledrelative}
\begin{equation}\label{eq:SGe}
    \begin{split}    
     \textup{SG}_{\mathcal{U}_e}(H) := &\left\{\rho_T(u,y)e^{\pm j\theta_T(u,y)} \mid \right.\\
    &\qquad\quad\left.u \in \mathcal{U}_e, \:y\in H(u), \:T>0\right\},
   \end{split}
\end{equation}
where $\rho_T(u,y) := \rho(P_Tu,P_Ty)$ and $\theta_T(u,y):=\theta(P_Tu,P_Ty)$, with $u,y\in \mathcal{L}_{2e}$, $T\geq 0$, and $\rho(\cdot,\cdot)$ and $\theta(\cdot,\cdot)$ are given in \eqref{eq:gain} and \eqref{eq:phase}, respectively. When $\mathcal{U}_e = \mathcal{L}_{2e}$ we write $\textup{SG}_e(H)$.

The connection between hard scaled graphs and hard IQCs, the latter of the form
\begin{align} \label{eq:hardIQC}
\int_{0}^{T}\begin{bmatrix} y(t) \\ u(t) \end{bmatrix}^\top \left(\Pi \otimes I_n\right)\begin{bmatrix} y(t) \\ u(t) \end{bmatrix} dt \geq 0, \quad \forall T\geq 0,
\end{align} 
can be established in a similar manner as the connection between soft scaled graphs and soft IQCs in Lemma~\ref{lem:connection_to_SG}. That is, a system $H:\mathcal{L}_{2e}\rightrightarrows\mathcal{L}_{2e}$ satisfies \eqref{eq:hardIQC} over all $u\in \mathcal{U}_e$ and $y\in H(u)$ if and only if $\textup{SG}_{\mathcal{U}_e}(H) \subset \mathcal{S}(\Pi)$, where $\mathcal{S}(\Pi)$ is given in \eqref{eq:S_PI_canonical}. The proof of this statement is equivalent to the proof of Lemma~\ref{lem:connection_to_SG}, in which $\|P_Tu\|$, $\|P_Ty\|$, and $\langle P_Tu,P_Ty\rangle$ are used. 

\subsection{Computations}
The main mechanism that is exploited in Theorems~\ref{th:LTI_approx}, \ref{th:improved_reset_KYP}, and \ref{th:improved_PWL_KYP} for obtaining over-approximations of soft scaled graphs, is the construction of a quadratic storage function $W(x) = x^\top Px$ that satisfies the dissipation inequality in \eqref{eq:def_dissipativity} under the posed LMI conditions, see also Section~\ref{sec:dissp}. For soft scaled graph computations, we did not pose constraints on sign-definiteness of $W$ (and, therefore, no constraints on $P$), since we exploited the assumption that $\lim_{t\to\infty}W(x(t))=0$. For hard scaled graph computation, the IQC in \eqref{eq:hardIQC} must hold for all $T\geq 0$. This can be realized by enforcing $W(x)\geq 0$ for all $x \in \mathbb{R}^m$, and thus enforce the right-hand side of the inequality in \eqref{eq:def_dissipativity} to be nonnegative for all $T\geq 0$. By adding the constraint $P\succeq 0$ to Theorems ~\ref{th:LTI_approx}, \ref{th:improved_reset_KYP}, and \ref{th:improved_PWL_KYP}, and the corresponding optimization problems, the hard scaled graph of a system $H$ is obtained via the same computational procedure as for soft scaled graphs. Note that in line with Remark~\ref{rem:PIpass}, matrices of the form in \eqref{eq:PIpasive} can be included for computations. This is important when working with, e.g., LTI integrators with transfer function $1/s$, which satisfy $\langle P_Tu,P_Ty\rangle \geq 0$, and thus their hard scaled graph covers the entire right-half complex plane. 

Summarizing: all foregoing connections between soft scaled graphs, IQCs and LMIs can be leveraged to hard scaled graphs by simply enforcing the extra constraint $P\succeq 0$.

\begin{figure*}[t]
	
	\centering
	\begin{subfigure}{0.32\linewidth}
    \centering
	\includegraphics[width=0.9\linewidth]{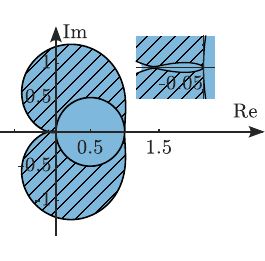}
	\caption{}
    \label{fig:LTI_hard}
	\end{subfigure}
	\begin{subfigure}{0.32\linewidth}
    \centering
	\includegraphics[width=0.9\linewidth]{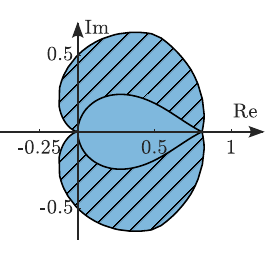}
	\caption{}
    \label{fig:reset_hard}
	\end{subfigure}
    \begin{subfigure}{0.32\linewidth}
    \centering
	\includegraphics[width=0.9\linewidth]{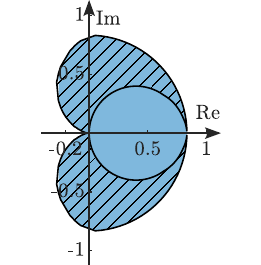}
	\caption{}
    \label{fig:_PWL_hard}
	\end{subfigure}
    
	\caption{Over-approximations of the hard scaled graphs in blue and the soft scaled graphs in the hatched black regions. (a) the LTI system from Example \ref{ex:LTI}, (b) the reset system from Example \ref{ex:reset} and (c) the PWL system from Example \ref{ex:PWL}.}
\label{fig:hard_example_sampled}
\end{figure*}

\subsection{Examples revisited}
To highlight the differences between soft and hard scaled graph computations, we reconsider Examples~\ref{ex:LTI}, \ref{ex:reset}, and~\ref{ex:PWL}. Over-approximations of the hard scaled graphs for each example are shown in Fig.~\ref{fig:hard_example_sampled} in blue. For comparison, the soft scaled graphs are shown by the hatched black regions. As expected, for $H=\left\{H_L,H_R,H_P\right\}$, $\textup{SG}(H) \subset \textup{SG}_e(H)$. It is interesting to note that the shape of $\textup{SG}(H)$ is equivalent to that of $\textup{SG}_e(H)$, and the only difference comes from the observation that the circular and ``teardrop'' inner shapes are included in the hard scaled graphs. Note that this also follows for the example in Remark~\ref{rem:1D}, in which the solution given by $P=-\frac{1}{2a}$ is excluded, leading to the interior of the disk with center $\lambda_c=\frac{1}{2a}$ and radius $r=\frac{1}{2a}$. For $a\to 0$, this region tends to the right-half plane, as expected.

\section{Conclusions}\label{sec:concl}
In this paper, we have presented a framework for computing (tight) over-approximations of scaled graphs for multivariable linear time-invariant, reset, and piecewise linear systems. Our approach exploits connections between LMIs, IQCs and scaled graphs, allowing to formulate the computation of (over-approximations of) so-called soft scaled graphs for these classes of systems as numerically tractable optimization procedures. For normal linear time-invariant systems, we show that the method yields the \emph{exact} scaled (relative) graph, thus introducing no conservatism in this case. Effectiveness of the method is shown through numerical experiments. We extended the framework to facilitate the computations of hard scaled graphs as well. For future work, we aim to extend the framework to broader classes of nonlinear systems, and leverage it to the design of nonlinear controllers, for instance via loop-shaping-inspired techniques \cite{chaffey_loop_2022} (see also Section~\ref{sec:inter}).

\bibliographystyle{plain}        
\bibliography{Sources}           

\end{document}